\documentclass[12pt]{article}
\usepackage{amssymb,amsmath,amsopn,amsfonts}
\usepackage[dvips]{color}
\usepackage{colordvi,multicol}
\usepackage{epsf}
\newfam\msbfam
\font\tenmsb=msbm10 \textfont\msbfam=\tenmsb \font\sevenmsb=msbm7
\scriptfont\msbfam=\sevenmsb \font\fivemsb=msbm5
\scriptscriptfont\msbfam=\fivemsb

\def\th#1{\vspace{1mm}\noindent{\bf #1}\quad}

\def\proof{\vspace{1mm}\noindent{\it Proof}\quad}

\numberwithin{equation}{section}

\def\bc{\begin{center}}
\def\ec{\end{center}}
\def\no{\noindent}
\def\hang{\hangindent\parindent}
\def\textindent#1{\indent\llap{\qquad #1\ \ \enspace}\ignorespaces}
\def\ref{\par\hang\textindent}

\begin{document}
\bigbreak

\title{ {\bf Ergodicity for the stochastic quantization problems on the 2D-torus
\thanks{Research supported in part  by NSFC (No.11301026, No.11401019)  and DFG through  CRC 701}\\} }
\author{{\bf Michael R\"{o}ckner}$^{\mbox{c},}$, {\bf Rongchan Zhu}$^{\mbox{a},\mbox{c}}$\thanks{
 Corresponding author} {\bf Xiangchan Zhu}$^{\mbox{b},\mbox{c}}$,
\date {}
\thanks{E-mail address:  roeckner@math.uni-bielefeld.de(M. R\"{o}ckner), zhurongchan@126.com(R. C. Zhu), zhuxiangchan@126.com(X. C. Zhu)}\\ \\
$^{\mbox{a}}$Department of Mathematics, Beijing Institute of Technology, Beijing 100081,
 China,\\
$^{\mbox{b}}$School of Science, Beijing Jiaotong University, Beijing 100044, China\\
$^{\mbox{c}}$ Department of Mathematics, University of Bielefeld, D-33615 Bielefeld, Germany,}

\maketitle

\noindent {\bf Abstract}

In this paper we study the  stochastic quantization problem on the two dimensional torus and establish  ergodicity for  the solutions. Furthermore, we prove a characterization of the $\Phi^4_2$ quantum field on the torus in terms of its density under translation. We  also deduce that the $\Phi^4_2$ quantum field on the torus is an extreme point in the set of all $L$-symmetrizing measures, where $L$ is the corresponding generator.

\vspace{1mm}
\no{\footnotesize{\bf 2000 Mathematics Subject Classification AMS}:\hspace{2mm} 60H15, 82C28}%Five letters(Mathematics Subject Classification 2000)
 \vspace{2mm}

\no{\footnotesize{\bf Keywords}:\hspace{2mm}  stochastic quantization problem, asymptotic coupling, space-time white noise, wick product}% At least three keywords

\section{Introduction}

In this paper we consider stochastic quantization equations on $\mathbb{T}^2$:
Let $H=L^2(\mathbb{T}^2)$:
\begin{equation}\aligned dX=&(AX-a_1:X^3:+a_2X)dt+dW(t),\\X(0)=&x,\endaligned\end{equation}
where $A:D(A)\subset H\rightarrow H$ is the linear operator
$$Ax=\Delta x-x,\quad D(A)=H^2(\mathbb{T}^2). $$
 $:x^3:$ means the renormalization of $x^3$ whose definition  will be  given in  Section 2, $a_1>0$ and $a_2$ is a real parameter. $W$ is the $L^2(\mathbb{T}^2)$-cylindrical $(\mathcal{F}_t)$-Wiener process defined on a probability space $(\Omega,\mathcal{F},P)$ equipped with a normal filtration $(\mathcal{F}_t)$.

This equation arises in the stochastic quantization of  Euclidean quantum field
theory. Consider the measure
$$\nu(d\phi)=ce^{-2\int:q(\phi):d\xi}\mu(d\phi),$$ where $q(\phi)=\frac{a_1}{4}\phi^4-\frac{a_2}{2}\phi^2$, $c$ is a normalization constant and $\mu$ is the Gaussian free field. The latter will be defined in Section 2 as well as the renormalization $:\cdot:$.  $\nu$ is called the $\Phi^4_2$-quantum field in Euclidean quantum field theory.  By heuristical calculations,  $\nu$ is an invariant measure for the solution to (1.1), which has been made rigorous in [RZZ15].  There have been many approaches to the problem of
giving a meaning to the above heuristic measure in the two dimensional case and the three dimensional case (see  [GRS75], [GlJ86] and the references
therein).
In [PW81]  Parisi and Wu proposed a program for Euclidean quantum field
theory of getting Gibbs states of classical
statistical mechanics as limiting distributions of stochastic processes,  especially those which are solutions to non-linear stochastic differential
equations. Then one can use the stochastic differential equations to study the properties of the Gibbs states. This
procedure is called stochastic field quantization (see [JLM85]). The $\Phi^4_2$ model is the simplest non-trivial Euclidean quantum field (see [GlJ86] and the reference therein). The issue of the stochastic quantization of the $\Phi^4_2$ model is to solve the equation (1.1) and to prove that the invariant measure is the limit of the time marginals as $t\rightarrow\infty$. The marginals converge to the Euclidean quantum field.

In [JLM85]  the existence of an ergodic, continuous, Markov
process having $\nu$ as an invariant measure has been proved, where $\nu$ is constructed above with $A$ changed to the Dirichlet Laplacian on a bounded domain. In fact, they consider the Markov process given by the solution to the following equation for $0<\varepsilon<\frac{1}{10}$
$$\aligned dX=&[-(-A)^\varepsilon X-(-A)^{-1+\varepsilon}(a_1:X^3:+a_2X)]dt+(-A)^{-\frac{1}{2}+\frac{\varepsilon}{2}}dW(t),\\X(0)=&x,\endaligned$$
which is easier to  solve than (1.1) (corresponding to the case $\varepsilon=1$) because of the regularization of the operator $A$. Moreover, they prove that the associated semigroup converges to $\nu$ in the $L^2$-sense.
 In [AR91]  weak solutions to (1.1) have been constructed by using the Dirichlet form approach in the finite  and infinite volume case.
 In [MR99] the stationary solution to (1.1) has also been considered in their general theory of martingale solutions for stochastic partial differential equations. In [DD03] again in the case of the torus, i.e. in finite volume,  Da Prato and Debussche define the Wick powers of solutions to the stochastic heat equation in the paths space and study a shifted equation instead of (1.1) in the finite volume case. They split the unknown $X$ into two parts:
$X=Y_1+Z_1$, where $Z_1(t)=\int_{-\infty}^te^{(t-s)A}dW(s)$. Observe that $Y_1$ is much smoother than $X$ and that in the stationary case
\begin{equation}:X^k:=\sum_{l=0}^kC_k^lY_1^l:Z_1^{k-l}:,\end{equation}
with $C_k^l=\frac{k!}{l!(k-l)!}$ and $:Z_1^{k-l}:$ being the Wick product,
which motivated them  to consider the following shifted equation:
\begin{equation}\aligned \frac{dY_1}{dt}=&AY_1-a_1\sum_{l=0}^{3} C_{3}^lY_1^l:Z_1^{3-l}:+a_2(Y_1+Z_1)\\Y_1(0)=&x-Z_1(0).\endaligned\end{equation}
They obtain  local existence and uniqueness of the solution $Y_1$ to (1.3) by a fixed point argument. By using the invariant measure $\nu$ they obtain a global solution to (1.1) by defining $X=Y_1+Z$ starting from almost every starting point. In [MW15] the authors consider the following equation  instead of (1.3):
 \begin{equation}\aligned \frac{dY}{dt}=&AY-a_1\sum_{l=0}^{3} C_{3}^lY^l:\bar{Z}^{3-l}:+a_2(Y+\bar{Z})\\Y(0)=&0,\endaligned\end{equation}
 where $\bar{Z}(t)=e^{tA}x+\int_0^te^{(t-s)A}dW(s)$ and $:\bar{Z}^{k-1-l}:$ will be defined later. We call (1.4) \emph{the shifted equation} for short.
 They obtain global existence and uniqueness of the solution to (1.4) directly from every starting point both in the finite and  infinite volume case.
Actually, (1.3) is equivalent to (1.4). For the solution $Y_1$ to (1.3), defining $Y(t):=Y_1(t)+e^{tA}Z_1(0)-e^{tA}x$, we can easily check that $Y$ is a solution to (1.4) by using the binomial formula (2.2) below.

In [RZZ15] we prove that  $X-\bar{Z}$, where $X$ is obtained by the Dirichlet form approach in [AR91] and $\bar{Z}(t)=\int_0^te^{(t-s)A}dW(s)+e^{tA}x$, also satisfies the shifted equation (1.4). Moreover, we obtain that the $\Phi^4_2$ quantum  field $\nu$ is an invariant measure for the process $X_0=Y+\bar{Z}$, where $Y$ is the unique solution to the shifted equation (1.4). It is natural to ask whether this invariant measure $\nu$ is the unique invariant measure for $X_0$. If $\nu$ is the unique invariant measure for $X_0$, then $\nu$ is the limiting distribution of the stochastic processes $X_0$.
 This problem is main point in the stochastic field quantization as we mentioned above in the $\Phi^4_2$ model on the torus.

 This problem has been studied in [AKR97] and the references therein. It is proved in [AKR97] that the stochastic quantization of a Guerra-Rosen-Simon Gibbs state on $\mathcal{S}'(\mathbb{R}^2)$ in infinite volume with polynomial interaction is ergodic if the Gibbs state is a pure phase. This result also holds for the finite volume case if one takes Dirichlet boundary conditions. Moreover, by [R86] we know that  $\nu$ constructed above with $A$ changed to a Dirichlet Laplacian on a bounded domain is a pure phase, which implies that the stochastic quantization of the Gibbs state is ergodic. However, the idea in [R86] and the results in [AKR97] cannot be applied for the torus. In this case we don't know whether $\nu$ is a pure phase. We also emphasize that it is not obvious that $\nu$ is a pure phase even if $\nu$ is absolutely continuous with respect to $\mu$. In this case, the zero set of $\frac{d\nu}{d\mu}$, i.e. $\{\frac{d\nu}{d\mu}=0\}$, which is hard to analyze analytically,  may divide the state space into different irreducible components, which immediately implies non-ergodicity, i.e. the existence of two invariant measures.   In this paper we study this problem using the techniques from SPDE. We analyze the shifted equation directly and obtain that $\nu$ is the unique invariant measure of $X_0$.

We also emphasize that  Dirichlet form theory is crucially used in [AKR97]. Hence for the Dirichlet boundary condition case, one can only obtain that the associated semigroup converges to the Gibbs state for quasi-every starting point. In our paper we analyze $X_0$ starting from every  point in $\mathcal{C}^\alpha$ for some $\alpha<0$, which will be defined in Section 2. As a result, we can conclude that the associated semigroup converges to $\nu$ for every starting point in $\mathcal{C}^\alpha$.

 \vskip.10in

\th{Theorem 1.1} $\nu$ is the unique  invariant probability measure for the process $X_0$.
Moreover, the associated semigroup $P_t$ converges to $\nu$ weakly in $\mathcal{C}^\alpha$, as $t$ goes to $\infty$.

 \vskip.10in

\th{Remark 1.2}As in [DD03], [RZZ15], one can replace the term $-:X^3:$ by any Wick polynomial
of odd degree with negative leading coefficient and obtain the same results in an analogous way.
 \vskip.10in

\th{Remark 1.3}  By [GlJ86] we know that for polynomials $\phi^4-\lambda\phi^2$ with $\lambda$ large enough,  the quantum fields in the infinite volume case may have different phases. We expect finding two different Gibbs state $\nu_1, \nu_2$ in this case such that they have similar property  as  in [GlJ86, Corollary 12.2.4]. If this is true, we can also obtain that these two states correspond to two different invariant measures for $X_0$ in the infinite volume case. However, so far  one only knows one state in the infinite volume case obtained in [GlJ86, Chapter 11] satisfying the property  in [GlJ86, Corollary 12.2.4].
 \vskip.10in

For the proof of Theorem 1.1, we use an argument from an abstract framework developed for
application to SDEs with delay [HMS11].
In general by applying a theorem in [HMS11] (see Theorem 4.1), we can reduce the problem of
uniqueness of the invariant measure to the convergence of solutions of (1.1) to solutions of an auxiliary system
 when time tends to infinity. However, in our case we cannot consider the solution to (1.1) obtained by Dirichlet form theory directly since it does not start from every point in some Polish space and the regularity of the solution to (1.1) is too rough to be controlled. Formally $:X^3:=X^3-\infty X$, which makes it more difficult to analyze the equation (1.1) directly.  Instead we consider the shifted equation (1.4) and do the required a-priori estimates for the solutions to (1.4). Correspondingly, we also construct an auxiliary system for the shifted form (see (3.3)). Moreover, to apply [HMS11] we have to construct a suitable set such that the generalized coupling has positive mass on it and the two solutions can converges to each other on this set when time tends to infinity .

 \vskip.10in
 As a consequence of Theorem 1.1 we can give a characterization of $\nu$ in terms of its density under translation:
  \vskip.10in

 \th{Theorem 1.4} $\nu$ is the unique probability measure such that the following hold

 i) $\nu$ is absolutely continuous with respect to $\mu$ with $\frac{d\nu}{d\mu}\in L^p(\mathcal{S}'(\mathbb{T}^2),\mu)$ for some $p>1$;

 ii) ("quasi-invariance") $\frac{d\nu(z+tk)}{d\nu(z)}=a_{tk}(z)=a_{tk}^0(z)a_{tk}^{a_1,a_2}(z)$ for $z\in H_2^{-1-\epsilon}, k\in C^\infty(\mathbb{T}^2)$, $t>0$
 with $$a_{tk}^0(z)=\exp[-t\langle (-\Delta +1)k,z\rangle-\frac{1}{2}t^2\langle (-\Delta +1)k,k\rangle]$$ and $$a_{tk}^{a_1,a_2}(z)=\exp[-\frac{a_1}{2}\sum_{i=0}^{3}C_4^it^{4-i}:z^i:(k^{4-i})+a_2\sum_{i=0}^{1}C_2^it^{2-i}:z^i:(k^{2-i})].$$
Here $H_2^{-1-\epsilon}$ for some $\epsilon>0$ is defined in  Section 2 and in the following $\langle \cdot,\cdot\rangle$ means the dualization between the elements in $C^\infty(\mathbb{T}^2)$ and $H_2^{-1-\epsilon}$, respectively. $:z^3:$ is a fixed version of the Wick power we define in Section 2.  By [AR91, Proposition 6.9] we can choose $z\rightarrow:z^3:$ as a measurable map from $H_2^{-1-\epsilon}$ to $H_2^{-1-\epsilon}$.
\vskip.10in
Similarly we obtain the following uniqueness result for the $L$-symmetrizing measures.
\vskip.10in
 \th{Theorem 1.5} $\nu$ is the unique probability measure such that the following hold:

 i) $\nu$ is absolutely continuous with respect to $\mu$ with $\frac{d\nu}{d\mu}\in L^p(\mathcal{S}'(\mathbb{T}^2),\mu)$ for some $p>1$;

 ii) $\int Lu v d\nu=\int Lv ud\nu$ for $u\in \mathcal{F}C_b^\infty$,
 where $$Lu(z)=\frac{1}{2}\textrm{Tr}(D^2u)(z)+\langle z,ADu\rangle-\langle a_1:z^3:-a_2z,Du(z)\rangle$$ for $z\in H_2^{-1-\epsilon}$ and $\mathcal{F}C_b^\infty$ is as defined in Section 4.

\vskip.10in
\th{Remark 1.6} i) From the proof of Theorems 1.4 and 1.5 we know that assuming  i) in Theorem 1.4 to hold, it follows that ii) in Theorems 1.4 and 1.5 are equivalent to the logarithmic derivative of $\nu$ along $k$ being given by
$$\beta_k=2\langle z,Ak\rangle-2\langle a_1:z^3:-a_2z,k\rangle,$$
for $z\in H_2^{-1-\epsilon}$, $k\in C^\infty(\mathbb{T}^2)$.
Here the  logarithmic derivative of a measure $\nu$ along $k$ is a $\nu$-integrable function $\beta_k$ such that the following integration by parts formula holds:
$$\int \frac{\partial u}{\partial k}d\nu=-\int\beta_k ud\nu.$$

\vskip.10in
Moreover, we can prove that $\nu$ is an extreme point of the following convex set.
\vskip.10in
\th{Corollary 1.7}  $\nu$ is an extreme point of the convex set $\mathcal{M}^a$, which denotes the set of all probability measures on $\mathcal{S}'(\mathbb{T}^2)$ satisfies (ii) in Theorem 1.4.
\vskip.10in

\th{Corollary 1.8}  $\nu$ is an extreme point of the convex set $\mathcal{G}$, which denotes the set of all probability measures on $\mathcal{S}'(\mathbb{T}^2)$ satisfies (ii) in Theorem 1.5.
\vskip.10in

\th{Remark 1.9} i) By [AKR97, Theorem 3.3] we know that $\nu$ being an extreme point of the convex set $\mathcal{M}^a$ is equivalent to $\nu$ being $C^\infty(\mathbb{T}^2)$-ergodic, which is also equivalent to the maximal Dirichlet form $(\mathcal{E}_\nu, D(\mathcal{E}_v))$ being irreducible.
For the definition of the maximal Dirichlet form $(\mathcal{E}_\nu, D(\mathcal{E}_v))$, we refer to [AKR97,  Section 3].

ii) Since the irreducibility is so crucial we recall here some characterizations of it in terms of the semigroup $(T_t)_{t>0}$ and generator $(L, D(L))$ of $(\mathcal{E}_\nu, D(\mathcal{E}_v))$. The following are equivalent:

1. $(\mathcal{E}_\nu, D(\mathcal{E}_v))$ is irreducible.

2. $(T_t)_{t>0}$ is irreducible, i.e., if $g\in L^2(\nu)$ such that $T_t(gf)=gT_tf$ for all $t>0$, $f\in L^2(\nu)$ then $g=const$.

3. If $g\in L^2(\nu)$ such that $T_tg=g$ for all $t>0$ then $g=const$.

4. $\int (T_tg-\int gd\nu)^2d\nu\rightarrow_{t\rightarrow\infty}0$ for all $g\in L^2(\nu)$.

5. If $u\in D(L)$ with $Lu=0$, then $u=const$.

Here we emphasize that we don't know whether the maximal Dirichlet form is the same as the minimal Dirichlet form defined in the proof of Theorem 1.4 below, which is the issue of the Markov uniqueness problem. If the maximal Dirichlet form is associated with a strong Markov process (i.e. is a quasi-regular Dirichlet form in the sense of [MR92]), then it is the same as the minimal Dirichlet form (see [RZZ15, Theorem 3.12]).

iii) The fact that the maximal Dirichlet form $(\mathcal{E}_\nu, D(\mathcal{E}_v))$ is irreducible can also be proved by a similar argument as in the proof of [BR95, Theorem 6.15] and by using Theorem 1.4.
\vskip.10in

We also want to mention that recently there
has arisen a renewed interest in SPDE’s related to such problems, particularly in
connection with Hairer's theory of regularity structures [Hai14] and related
work by Imkeller, Gubinelli, Perkowski in [GIP13]. By using these theories one can obtain local existence and uniqueness of solutions to (1.1) in the three dimensional case (see [Hai14, CC13]).

This paper is organized as follows: In Section 2 we collect some results related to Besov spaces and we recall some basic facts on   Wick powers. In Section 3 we prove the necessary a-priori estimates of solutions to (1.1). In Section 4 we prove Theorems 1.1, 1.4 and 1.5.

\section{Preliminaries}
\subsection{Notations and some useful estimates}

In the following we recall the definitions  of Besov spaces. For a general introduction to the theory we refer to [BCD11, Tri78, Tri06].
 The space of real valued infinitely differentiable functions of compact support is denoted by $\mathcal{D}(\mathbb{R}^d)$ or $\mathcal{D}$. The space of Schwartz functions is denoted by $\mathcal{S}(\mathbb{R}^d)$. Its dual, the space of tempered distributions is denoted by $\mathcal{S}'(\mathbb{R}^d)$. The Fourier transform and the inverse Fourier transform are denoted by $\mathcal{F}$ and $\mathcal{F}^{-1}$.

 Let $\chi,\theta\in \mathcal{D}$ be nonnegative radial functions on $\mathbb{R}^d$, such that

i. the support of $\chi$ is contained in a ball and the support of $\theta$ is contained in an annulus;

ii. $\chi(\xi)+\sum_{j\geq0}\theta(2^{-j}\xi)=1$ for all $\xi\in \mathbb{R}^d$.

iii. $\textrm{supp}(\chi)\cap \textrm{supp}(\theta(2^{-j}\cdot))=\emptyset$ for $j\geq1$ and $\textrm{supp}\theta(2^{-i}\cdot)\cap \textrm{supp}\theta(2^{-j}\cdot)=\emptyset$ for $|i-j|>1$.

We call such a pair $(\chi,\theta)$ dyadic partition of unity, and for the existence of dyadic partitions of unity see [BCD11, Proposition 2.10]. The Littlewood-Paley blocks are now defined as
$$\Delta_{-1}u=\mathcal{F}^{-1}(\chi\mathcal{F}u)\quad \Delta_{j}u=\mathcal{F}^{-1}(\theta(2^{-j}\cdot)\mathcal{F}u).$$

For $\alpha\in\mathbb{R}$, $p,q\in [1,\infty]$, $u\in\mathcal{D}$ we define
$$\|u\|_{B^\alpha_{p,q}}:=(\sum_{j\geq-1}(2^{j\alpha}\|\Delta_ju\|_{L^p})^q)^{1/q},$$
with the usual interpretation as $l^\infty$ norm in case $q=\infty$. The Besov space $B^\alpha_{p,q}$ consists of the completion of $\mathcal{D}$ with respect to this norm and the H\"{o}lder-Besov space $\mathcal{C}^\alpha$ is given by $\mathcal{C}^\alpha(\mathbb{R}^d)=B^\alpha_{\infty,\infty}(\mathbb{R}^d)$. For $p,q\in [1,\infty)$,
$$B^\alpha_{p,q}(\mathbb{R}^d)=\{u\in\mathcal{S}'(\mathbb{R}^d):\|u\|_{B^\alpha_{p,q}}<\infty\}.$$
$$\mathcal{C}^\alpha(\mathbb{R}^d)\varsubsetneq \{u\in\mathcal{S}'(\mathbb{R}^d):\|u\|_{\mathcal{C}^\alpha(\mathbb{R}^d)}<\infty\}.$$
We point out that everything above and everything that follows can be applied to distributions on the torus (see [S85, SW71]). More precisely, let $\mathcal{S}'(\mathbb{T}^d)$ be the space of distributions on $\mathbb{T}^d$.  Besov spaces on the torus with general indices $p,q\in[1,\infty]$ are defined as
the completion of $C^\infty(\mathbb{T}^2)$ with respect to the norm $$\|u\|_{B^\alpha_{p,q}(\mathbb{T}^d)}:=(\sum_{j\geq-1}(2^{j\alpha}\|\Delta_ju\|_{L^p(\mathbb{T}^d)})^q)^{1/q},$$
and the H\"{o}lder-Besov space $\mathcal{C}^\alpha$ is given by $\mathcal{C}^\alpha=B^\alpha_{\infty,\infty}(\mathbb{T}^d)$.  We write $\|\cdot\|_{\alpha}$ instead of $\|\cdot\|_{B^\alpha_{\infty,\infty}(\mathbb{T}^d)}$ in the following for simplicity.  For $p,q\in[1,\infty)$
$$B^\alpha_{p,q}(\mathbb{T}^d)=\{u\in\mathcal{S}'(\mathbb{T}^d):\|u\|_{B^\alpha_{p,q}(\mathbb{T}^d)}<\infty\}.$$
\begin{equation}\mathcal{C}^\alpha\varsubsetneq \{u\in\mathcal{S}'(\mathbb{T}^d):\|u\|_{\alpha}<\infty\}.\end{equation}

 \vskip.10in

In this part we give estimates on the torus for later use.
 Set $\Lambda= (-A)^{\frac{1}{2}}$. For $s\geq0, p\in [1,+\infty]$ we use $H^{s}_p$ to denote the subspace of $L^p(\mathbb{T}^d)$, consisting of all  $f$   which can be written in the form $f=\Lambda^{-s}g, g\in L^p(\mathbb{T}^d)$ and the $H^{s}_p$ norm of $f$ is defined to be the $L^p$ norm of $g$, i.e. $\|f\|_{H^{s}_p}:=\|\Lambda^s f\|_{L^p(\mathbb{T}^d)}$.

\vskip.10in
To study (1.1) in the finite volume case, we  will need several important properties of Besov spaces on the torus and we recall the following Besov embedding theorems on the torus first (c.f. [Tri78, Theorem 4.6.1], [GIP13, Lemma 41]):
\vskip.10in
 \th{Lemma 2.1} (i) Let $1\leq p_1\leq p_2\leq\infty$ and $1\leq q_1\leq q_2\leq\infty$, and let $\alpha\in\mathbb{R}$. Then $B^\alpha_{p_1,q_1}(\mathbb{T}^d)$ is continuously embedded in $B^{\alpha-d(1/p_1-1/p_2)}_{p_2,q_2}(\mathbb{T}^d)$.

 (ii) Let $s\geq0$, $1<p<\infty$, $\epsilon>0$. Then
 $ H^{s+\epsilon}_p\subset B^{s}_{p,1}(\mathbb{T}^d)\subset B^{s}_{1,1}(\mathbb{T}^d)$.

 (iii) Let $1\leq p_1\leq p_2<\infty$ and let $\alpha\in\mathbb{R}$. Then $H^\alpha_{p_1}$ is continuously embedded in $H^{\alpha-d(1/p_1-1/p_2)}_{p_2}$.

Here  $\subset$ means that the embedding is continuous and dense.

\vskip.10in

We recall the following Schauder estimates, i.e. the smoothing effect of the heat flow, for later use.

\vskip.10in
\th{Lemma 2.2}([GIP13, Lemma 47]) (i) Let $u\in B^{\alpha}_{p,q}(\mathbb{T}^d)$ for some $\alpha\in \mathbb{R}, p,q\in [1,\infty]$. Then for every $\delta\geq0$
$$\|e^{tA}u\|_{B^{\alpha+\delta}_{p,q}(\mathbb{T}^d)}\lesssim t^{-\delta/2}\|u\|_{B^{\alpha}_{p,q}(\mathbb{T}^d)},$$
where the constant we omit is independent of $t$.

(ii) Let $\alpha\leq \beta\in\mathbb{R}$. Then
$$\|(1-e^{tA})u\|_{\alpha}\lesssim t^{\frac{\beta-\alpha}{2}}\|u\|_{\beta}.$$

\vskip.10in
One can  extend the multiplication on suitable Besov spaces and also  have the duality properties of Besov spaces from [Tri78, Chapter 4]:
\vskip.10in

\th{Lemma 2.3} (i) Let $\alpha,\beta\in\mathbb{R}$ and $p, p_1, p_2, q\in [1,\infty]$ be such that $$\frac{1}{p}=\frac{1}{p_1}+\frac{1}{p_2}.$$
The bilinear map $(u; v)\mapsto uv$
extends to a continuous map from ${B}^\alpha_{p_1,q}\times {B}^\beta_{p_2,q}$ to ${B}^{\alpha\wedge\beta}_{p,q}$  if $\alpha+\beta>0$.

(ii) Let $\alpha\in (0,1)$, $p,q\in[1,\infty]$, $p'$ and $q'$ be their conjugate exponents, respectively. Then the mapping  $(u; v)\mapsto \int uvdx$  extends to a continuous bilinear form on $B^\alpha_{p,q}(\mathbb{T}^d)\times B^{-\alpha}_{p',q'}(\mathbb{T}^d)$.

\vskip.10in
We recall the following interpolation inequality and  multiplicative inequality for the elements in $H^s_p$ (cf. [Tri78, Theorem 4.3.1], [Re95, Lemma A.4], [RZZ15a, Lemma 2.1]):
 \vskip.10in

\th{Lemma 2.4} (i)  Suppose that $s\in (0,1)$ and $p\in (1,\infty)$. Then for $u\in H^1_p$
$$\|u\|_{H^s_p}\lesssim \|u\|_{L^p(\mathbb{T}^d)}^{1-s}\|u\|_{H^1_p}^s.$$

(ii) Suppose that $s>0$ and $p\in (1,\infty)$. If $u,v\in C^{\infty}(\mathbb{T}^2)$ then
$$\|\Lambda^s(uv)\|_{L^p(\mathbb{T}^d)}\lesssim\|u\|_{L^{p_1}(\mathbb{T}^d)}\|\Lambda^sv\|_{L^{p_2}(\mathbb{T}^d)}+\|v\|_{L^{p_3}(\mathbb{T}^d)}
\|\Lambda^su\|_{L^{p_4}(\mathbb{T}^d)},$$
with $p_i\in (1,\infty], i=1,...,4$ such that
$$\frac{1}{p}=\frac{1}{p_1}+\frac{1}{p_2}=\frac{1}{p_3}+\frac{1}{p_4}.$$
\vskip.10in

\subsection{Wick power}
In the following we recall  the definition of Wick powers. Let $\mu=N(0,\frac{1}{2}(-\Delta+1)^{-1}):=N(0,C)$.
 \vskip.10in
\textbf{Wick power on $L^2(\mathcal{S}'(\mathbb{T}^2),\mu)$}

In fact $\mu$ is a measure supported on $\mathcal{S}'(\mathbb{T}^2)$. We have the well-known (Wiener-It\^{o}) chaos decomposition $$L^2(\mathcal{S}'(\mathbb{T}^2),\mu)=\bigoplus_{n\geq0}\mathcal{H}_n.$$
Now we define the Wick power by using approximations: for $\phi\in \mathcal{S}'(\mathbb{T}^2)$ define $$\phi_\varepsilon:=\rho_\varepsilon*\phi,$$ with $\rho_\varepsilon$ an approximate delta function,
$$\rho_\varepsilon(\xi)=\varepsilon^{-2}\rho(\frac{\xi}{\varepsilon})\in \mathcal{D}, \int \rho=1.$$
Here the convolution means that we view $\phi$ as a periodic distribution in $\mathcal{S}'(\mathbb{R}^2)$.
 For every $n\in\mathbb{N}$ we set $$:\phi_\varepsilon^n:_C:=c_\varepsilon^{n/2}P_n(c_\varepsilon^{-1/2}\phi_\varepsilon),$$
where $P_n,n=0,1,...,$ are the Hermite polynomials defined by the formula
$$P_n(x)=\sum_{j=0}^{[n/2]}(-1)^j\frac{n!}{(n-2j)!j!2^j}x^{n-2j},$$
and
$c_\varepsilon:=\int\phi^2_\varepsilon\mu(d\phi)=\int\int \bar{G}(\xi_1-\xi_2)\rho_\varepsilon(\xi_2)d\xi_2\rho_\varepsilon(\xi_1)d\xi_1
=\|\bar{K}_\varepsilon\|_{L^2(\mathbb{R}\times \mathbb{T}^2 )}^2$. Then $:\phi_\varepsilon^n:_C\in \mathcal{H}_n$.
Here and in the following $\bar{G}$ is the Green function associated with $-A$ on $\mathbb{T}^2$ and let $\bar{K}(t,\xi)$ be such that $\bar{K}(t,\xi_1-\xi_2)$ is the heat kernel associated with $A$ on $\mathbb{T}^2$ and $\bar{K}_\varepsilon=\bar{K}*\rho_\varepsilon$ with $*$ means  convolution in space and we view $\bar{K}$ as a periodic function on $\mathbb{R}^2$.

For Hermite polynomial $P_n$ we have for $s,t\in\mathbb{R}$
\begin{equation}P_n(s+t)=\sum_{m=0}^nC_n^mP_m(s)t^{n-m},\end{equation}
where $C_n^m=\frac{n!}{m!(n-m)!}$.

 \vskip.10in
A direct calculation yields the following:
 \vskip.10in

\th{Lemma 2.5}([RZZ15, Lemma 3.1]) Let $\alpha<0$, $n\in\mathbb{N}$ and $p>1$. $:\phi_\varepsilon^n:_C$ converges to some element
in $L^p(\mathcal{S}'(\mathbb{T}^2),\mu;\mathcal{C}^{\alpha})$. This limit is called $n$-th Wick power of $\phi$ with respect to the covariance $C$ and denoted by $:\phi^n:_C$.

\vskip.10in
Now we introduce the following probability measure.   Let $$\nu=c\exp{(-\frac{1}{2}\int_{\mathbb{T}^2}(a_1:\phi^4:_C-2a_2:\phi^2:_C)d\xi)}\mu,$$
where $c$ is a normalization constant. Then by  [GlJ86, Sect. 8.6] for every $p\in [1,\infty)$, $\varphi(\phi):=\exp{(-\frac{1}{2}\int_{\mathbb{T}^2}(a_1:\phi^4:_C-2a_2:\phi^2:_C)d\xi)}\in L^p(\mathcal{S}'(\mathbb{T}^2),\mu)$.

 \vskip.10in
\textbf{Wick power on a fixed probability space}

Now we fix a probability space $(\Omega,\mathcal{F},P)$ equipped with a normal filtration $(\mathcal{F}_t)$ and $W$ is an $L^2(\mathbb{T}^2)$-cylindrical
$(\mathcal{F}_t)$-Wiener process. We also have the well-known (Wiener-It\^{o}) chaos decomposition
$$L^2(\Omega,\mathcal{F},P)=\bigoplus_{n\geq0}\mathcal{H}'_n.$$In the following we set $Z(t)=\int_0^t e^{(t-s)A}dW(s)$, and we can also define Wick powers of
$Z(t)$  by approximations: Let $Z_\varepsilon(t,\xi)=\int_0^t \langle\bar{K}_\varepsilon(t-s,\xi-\cdot),dW(s)\rangle$. Here $\langle\cdot,\cdot\rangle$ means  inner product in $L^2(\mathbb{T}^2)$.
 \vskip.10in

\th{Lemma 2.6}([RZZ15, Lemma 3.4]) For  $\alpha<0$, $n\in\mathbb{N}$ and $t>0$, $:Z_\varepsilon^n(t)::=c_\varepsilon^{\frac{n}{2}}P_n(c_\varepsilon^{-\frac{1}{2}}Z_\varepsilon (t))$ converges  in $L^p(\Omega,C((0,T];\mathcal{C}^{\alpha}))$. Here the norm for $C((0,T];\mathcal{C}^{\alpha})$ is $\sup_{t\in[0,T]}t^\delta\|\cdot\|_\alpha$ for  $\delta>0$. The limit is called Wick power of $Z(t)$ with respect to the covariance $C$ and denoted by $:Z^n(t):$.

\vskip.10in

Now following the technique in [MW15] we combine the initial value part  with the Wick powers by using (2.2).
We set $V(t)=e^{tA}x$, $x\in{\mathcal{C}}^{\alpha}$ for $\alpha<0$ and
$$\bar{Z}_x(t)=Z(t)+V(t),$$
and for $n=2,3$,
$$:\bar{Z}^n_x(t):=\sum_{k=0}^nC_n^kV(t)^{n-k}:Z^k(t):.$$
By Lemma 2.2 we know that $V\in C([0,T],\mathcal{C}^\alpha)$  and  $V\in C((0,T],\mathcal{C}^\beta)$ for $\beta>-\alpha$ with the norm $\sup_{t\in[0,T]}t^{\frac{\beta-\alpha}{2}}\|\cdot\|_{\beta}$. Moreover,
\begin{equation}t^{\frac{\beta-\alpha}{2}}\|V(t)\|_{\beta}\lesssim \|x\|_{\alpha},\end{equation}
for $\beta>-\alpha$. Then by [RZZ15, Lemmas 3.5] we have $\bar{Z}_x\in L^p(C((0,T],\mathcal{C}^\alpha))$.
\vskip.10in
By (2.3) and Lemma 2.3 it is easy to obtain the following result:
\vskip.10in

\th{Lemma 2.7} Let $\alpha<0, x\in\mathcal{C}^\alpha$. Then  we have for $t>0$ and any $\varepsilon>0$
$$\aligned\|\bar{Z}_x(t)\|_{\alpha}\leq& \|Z(t)\|_{\alpha}+\|x\|_{\alpha},\\\|:\bar{Z}_x^2(t):\|_{\alpha}\leq& C[ \|:Z^2(t):\|_{\alpha}+t^{\alpha-\varepsilon}\|x\|_{\alpha}\|Z(t)\|_{\alpha}+t^{\alpha-\varepsilon}\|x\|_{\alpha}^2],
\\\|:\bar{Z}_x^3(t):\|_{\alpha}\leq& C[ \|:Z^3(t):\|_{\alpha}+t^{2\alpha-\varepsilon}\|x\|_{\alpha}^2\|Z(t)\|_{\alpha}
\\&+t^{\alpha-\varepsilon}\|x\|_{\alpha}\|:Z^2(t):\|_{\alpha}+t^{2\alpha-\varepsilon}\|x\|_{\alpha}^3].\endaligned$$

\section{The necessary a-priori estimates}

Now we follow [MW15, RZZ15] and give the existence and uniqueness of solutions to (1.1) by considering the shifted equation.

We fix  $\alpha<0$ with $-\alpha$ small enough, in the following. For $\underline{{Z}}=({Z},:{Z}^2:,:{Z}^3: )$, $0<\delta<-\alpha$,
define $$\|\underline{{Z}}\|_{\mathfrak{L}_T}:=\sup_{0\leq t\leq T}\big(\|Z(t)\|_\alpha, t^{\delta}\|:Z^2(t):\|_\alpha, t^{\delta}\|:Z^3(t):\|_\alpha  \big)$$
and $$\Omega_0=\{Z\in C([0,T];{\mathcal{C}}^{\alpha}), :Z^n:\in C((0,T]; {\mathcal{C}}^{\alpha}) ,\|\underline{{Z}}\|_{\mathfrak{L}_T}<\infty \textrm{ for } n=2,3, T>0\}.$$
Then  $$P[\Omega_0]=1.$$

We also introduce the following notations: for $Y\in C([0,T],\mathcal{C}^\beta),
Z\in C([0,T];{\mathcal{C}}^{\alpha})$ with $\beta>-\alpha$
$$\Psi(Y,Z):=-a_1(3Y^2Z+3Y:Z^2:+:Z^3:)+a_2(Y+Z).$$
By [RZZ15, Theorem 3.9] we know that the solution to (1.1) can be written as $X=Y+e^{tA}x+Z=Y+\bar{Z}_x$ for $x\in\mathcal{C}^\alpha$,
where $Y$  satisfies the following shift  equation:
\begin{equation}dY=[AY-a_1Y^3+\Psi(Y,{\bar{Z}_x})]dt,\quad Y(0)=0.\end{equation}
Here and in the following (3.1) and other equations are interpreted in the mild sense:
$$Y(t)=\int_0^te^{(t-s)A}[-a_1Y^3+\Psi(Y,{\bar{Z}_x})](s)ds.$$
As a result, $X$ is a mild solution to the following equation
\begin{equation}dX=[AX-a_1(X-\bar{Z}_{x})^3+\Psi(X-\bar{Z}_x,{\bar{Z}_x})]dt+dW,\quad X(0)=x.\end{equation}
That is to say:
\vskip.10in

\th{Theorem 3.1} For $\omega\in \Omega_0, x\in \mathcal{C}^\alpha$, there exists exactly one mild solution $X(\omega)\in C([0,\infty);\mathcal{C}^\alpha)$ to
the equation (3.2) satisfying $(X-\bar{Z}_x)(\omega)\in C([0,\infty);\mathcal{C}^\beta)$ for some $\beta>-\alpha>0$.

\proof  For $\omega\in\Omega_0$,  $x\in \mathcal{C}^\alpha$,
by Lemma 2.7 we know that for $n=2,3$,
$$\bar{Z}_x\in C([0,\infty); {\mathcal{C}}^{\alpha}),\quad :\bar{Z}_x^n:\in C((0,\infty); {\mathcal{C}}^{\alpha}),$$
and for every $T>0$
$$\sup_{0\leq t\leq T}\big[\|\bar{Z}_x(t)\|_\alpha, t^{-2\alpha+\delta}\|:\bar{Z}_x^2(t):\|_\alpha, t^{-2\alpha+\delta}\|:\bar{Z}_x^3(t):\|_\alpha  \big]<\infty.$$
Then  [MW15, Theorem 6.5] implies that for $\omega \in \Omega_0$ there exist exactly one solution $Y(\omega)\in C([0,\infty);\mathcal{C}^{\beta})$ for some $\beta>-\alpha$ satisfying (3.2) in the mild sense. From this we can conclude the result easily.  $\hfill\Box$

\vskip.10in
\th{Remark 3.2} Here for $\omega\in\Omega_0$ $X$ is an $\omega$-wise mild solution to the equation (3.2),
whose definition is  stronger than the usual (probabilistically) mild solution to the stochastic differential equation (3.2).

\vskip.10in

Now we can define the semigroup associated with $X$ and obtain an invariant measure for $X$:
 For $t\geq0, f\in\mathcal{B}_b(\mathcal{C}^\alpha)$, define $$P_tf(x):=Ef(X(t,x))$$ for the solution $X(t,x)$ to (3.2)
 with initial value  $x\in\mathcal{C}^\alpha$. By Theorem 3.1 we obtain that $X$ is a Markov process with $(P_t)_{t\geq0}$ as the associated semigroup.
 By [RZZ15, Theorem 3.10] we know that $\nu$ is an invariant measure for $X$. In the next section we will prove that $\nu$ is the unique invariant measure for $X$. To prove this,
 we introduce the following equation, which has a new dissipation term compared to (3.2).
\vskip.10in

For  given $x_0, x_1\in \mathcal{C}^{\alpha}$ consider the following equation
\begin{equation}\frac{d}{dt}\tilde{Y}=A\tilde{Y}-\lambda(\tilde{Y}-Y+e^{tA}(x_1-x_0))-a_1\tilde{Y}^3+\Psi(\tilde{Y},{\bar{Z}_{x_1}}),\quad \tilde{Y}(0)=0,\end{equation}
where $\lambda>1$ will be determined later.

By similar arguments as the proof of [MW15, Theorem 6.5] and [RZZ15, Theorem 3.10] we can easily derive the following result:
\vskip.10in

\th{Theorem 3.3} For  $x_0, x_1\in\mathcal{C}^\alpha$ and $\omega\in \Omega_0$, there exists a unique  mild solution $\tilde{Y}\in C([0,\infty);\mathcal{C}^\beta)$ to
the equation (3.3).
\vskip.10in

Define $\tilde{X}:=\tilde{Y}+e^{tA}x_1+Z=\tilde{Y}+\bar{Z}_{x_1}$ with $\tilde{Y}$ obtained in Theorem 3.3. Then $\tilde{X}$ satisfies the following equation in the mild sense:
\begin{equation}d\tilde{X}=[A\tilde{X}-\lambda(\tilde{X}-X)-a_1(\tilde{X}-\bar{Z}_{x_1})^3+\Psi(\tilde{X}-\bar{Z}_{x_1},{\bar{Z}_{x_1}})]dt+dW,\quad
\tilde{X}(0)=x_1.\end{equation}

Now we give the necessary a-priori estimates for the solutions to (3.1) and (3.3) for  later use. We will derive  $L^p$-norm estimates for the solutions to (3.1) and (3.3) respectively.  We can get the $\|\cdot\|_{L^p}$-norm estimate directly, but only with $p$ depending on $\alpha$, which is not enough for later use. Hence we first obtain the $\|\cdot\|_{L^p}$-norm estimates from $1$ to $t$ for every $p>1$. Then we estimate $\|Y(1)\|_{\beta}$ for $\beta>0$ by the following Steps 2 and 3.
\vskip.10in

\th{Lemma 3.4} Suppose that $Y$ is the solution to (3.1) with $x=x_0$. For every even $p>1$,  there exist constants $C(p), C(\|\underline{Z}\|_{\mathfrak{L}_1}, \|x_0\|_{\alpha})>0,\gamma(p)>1$ independent of $\omega, t$ such that for every $t\geq1$ and $\omega\in\Omega_0$
$$\aligned&\|Y(t)\|_{L^p}^p+\int_1^t\|Y(s)\|_{L^p}^pds+\int_1^t\|Y^{p-2}|\nabla Y|^2(s)\|_{L^1}ds\\\leq& C(\|\underline{Z}\|_{\mathfrak{L}_1}, \|x_0\|_{\alpha})+ C(p)\int_1^t (1+\|x_0\|_\alpha^{\gamma(p)}+\|{Z}\|_{\alpha}^{\gamma(p)}+\sum_{n=2}^3\|:{Z}^n:\|_{\alpha}^{\gamma(p)})ds.\endaligned$$
Here $C(\|\underline{Z}\|_{\mathfrak{L}_1}, \|x_0\|_{\alpha})$ means a constant depending on $\|\underline{Z}\|_{\mathfrak{L}_1}, \|x_0\|_{\alpha}$.

\proof We write the proof for $Y$ and $Z$ directly which is a bit informal,  but it  can be made rigorous by replacing $Z$ by $Z_\varepsilon$ and taking the limit  as in the proof of [RZZ15, Theorem 6.5].

\textbf{Step 1} We first prove that for every even $p>1$, $t\geq1$ there exists $\gamma(p)>2$ such that
\begin{equation}\aligned&\|Y(t)\|_{L^p}^p+\int_{1}^t\|Y^{p-2}|\nabla Y|^2\|_{L^1}^pds+\int_{1}^t\|Y(s)\|_{L^p}^pds
\\\leq& \|Y(1)\|_{L^p}^p+ C(p)\int_{1}^t (1+\|x_0\|_\alpha^{\gamma(p)}+\|Z\|_{\alpha}^{\gamma(p)}+\sum_{n=2}^3\|:Z^n:\|_{\alpha}^{{\gamma(p)}})ds.
\endaligned\end{equation}

Testing against $Y^{p-1}$, we  have that for $t\geq1$, even $p>1$,
$$\aligned&\frac{1}{p}\|Y(t)\|_{L^p}^p+\int_1^t[(p-1)\langle\nabla Y(s),
Y(s)^{p-2}\nabla Y(s)\rangle+a_1\|Y(s)^{p+2}\|_{L^1}]ds\\&=-\int_1^t[\|Y(s)\|_{L^p}^p+\langle \Psi(Y(s),{\bar{Z}_{x_0}}(s)),Y(s)^{p-1}\rangle]ds+\frac{1}{p}\|Y(1)\|_{L^p}^p.\endaligned$$
Now we have $\langle\nabla Y(s), Y(s)^{p-2}\nabla Y(s)\rangle$ and $\|Y(s)^{p+2}\|_{L^1}$ on the left hand side of the equality, which can be used to control the right hand side of the above equation.  By similar calculations as in the proof of  [MW15, Theorem 6.4] and [RZZ15, Theorem 3.10] we  deduce that there exists $ \gamma_0>1$ such that
$$\aligned&|\langle \Psi(Y(s),{\bar{Z}_{x_0}}(s)), Y(s)^{p-1}\rangle|\\\leq & C(p) (1+\|\bar{Z}_{x_0}\|_{\alpha}^{\gamma_0}+\sum_{n=2}^3\|:\bar{Z}^n_{x_0}:\|_{\alpha}^{p+2})
+\frac{1}{2} (a_1\|Y\|_{L^{p+2}}^{p+2}+\| Y^{p-2}|\nabla Y|^2\|_{L^{1}}),\endaligned$$
which implies that
\begin{equation}\aligned&\frac{1}{p}\|Y(t)\|_{L^p}^p
+\frac{1}{2}\int_{1}^t[(p-1)\langle\nabla Y(s), Y(s)^{p-2}\nabla Y(s)\rangle+a_1\|Y(s)^{p+2}\|_{L^1}]ds
\\\leq&C(p)\int_{1}^t (1+\|\bar{Z}_{x_0}\|_{\alpha}^{\gamma_0}+\sum_{n=2}^3\|:\bar{Z}^n_{x_0}:\|_{\alpha}^{p+2})ds+\frac{1}{p}\|Y(1)\|_{L^p}^p\\\leq&C(p)\int_{1}^t (1+\|x_0\|_\alpha^{\gamma(p)}+\|Z\|_{\alpha}^{\gamma(p)}+\sum_{n=2}^3\|:Z^n:\|_{\alpha}^{\gamma(p)})ds+\frac{1}{p}\|Y(1)\|_{L^p}^p.\endaligned\end{equation}
Here $\gamma(p)=3(p+2)\vee \gamma_0$ and we used Lemma 2.7 in the last inequality.
Now (3.5) follows.

\textbf{Step 2} We prove that for even $p>1$, with $-2\alpha(p+2)<1$,
\begin{equation}\sup_{0\leq t\leq 1}\|Y(t)\|_{L^p}^p+\int_0^1\|Y^{p-2}|\nabla Y|^2\|_{L^1}ds\leq C(p, \|\underline{Z}\|_{\mathfrak{L}_1}, \|x_0\|_{\alpha}).\end{equation}
By similar arguments as in Step 1 we have for $0\leq t\leq 1$
$$\aligned&\frac{1}{p}\|Y(t)\|_{L^p}^p+\frac{1}{2}\int_{0}^t[(p-1)\langle\nabla Y(s), Y(s)^{p-2}\nabla Y(s)\rangle
+a_1\|Y(s)^{p+2}\|_{L^1}]ds\\\leq&C(p)\int_{0}^t (1+\|\bar{Z}_{x_0}\|_{\alpha}^{\gamma_0}+\sum_{n=2}^3\|:\bar{Z}^n_{x_0}:\|_{\alpha}^{p+2})ds\\\leq& C(p, \|\underline{Z}\|_{\mathfrak{L}_1}, \|x_0\|_{\alpha}).\endaligned$$
Here we used Lemma 2.7 and $-2\alpha(p+2)<1$ in the last inequality.
Now (3.7) follows.

\textbf{Step 3} We prove that for $0<\beta<\frac{1}{2}+\alpha$
\begin{equation}\|Y(1)\|_{\beta}\leq C(\|\underline{Z}\|_{\mathfrak{L}_1}, \|x_0\|_{\alpha}).\end{equation}
Since $Y$ satisfies the mild equation, we have
 $$\aligned\|Y(1)\|_{\beta}\leq&C\int_0^1(1-s)^{-\frac{\beta}{2}-\frac{1}{2}}\|Y^3\|_{L^2}ds
+C\int_0^1(1-s)^{-\frac{\beta+1/2-\alpha}{2}}[\|Y^2\bar{Z}_{x_0}\|_{B^\alpha_{4,\infty}}+\|Y:\bar{Z}_{x_0}^2:\|_{B^\alpha_{4,\infty}}]ds
\\&+C\int_0^1(1-s)^{-(\beta-\alpha)/2}\|:\bar{Z}_{x_0}^3:\|_\alpha ds,
\endaligned$$
where  we used Lemmas 2.1, 2.2 to deduce $\|e^{tA}x\|_{\beta}\leq Ct^{-\frac{\beta}{2}-\frac{1}{2}}\|x\|_{L^2}$ and $\|\cdot\|_{\alpha-\frac{1}{2}}
\leq C\|\cdot\|_{B^\alpha_{4,\infty}}$. For the first term on the right hand side, we can use (3.7) for $p=6$ to control it by $C(\|\underline{Z}\|_{\mathfrak{L}_1}, \|x_0\|_{\alpha})$. Using Lemma 2.7 we can control the third term  by $C(\|\underline{Z}\|_{\mathfrak{L}_1}, \|x_0\|_{\alpha})$. Now we come to the second term:
$$\aligned &\int_0^1(1-s)^{-\frac{\beta+1/2-\alpha}{2}}[\|Y^2\bar{Z}_{x_0}\|_{B^\alpha_{4,\infty}}+\|Y:\bar{Z}_{x_0}^2:\|_{B^\alpha_{4,\infty}}]ds\\\leq& C\int_0^1(1-s)^{-\frac{\beta+1/2-\alpha}{2}}[\|Y^2\|_{B^\beta_{4,\infty}}
\|\bar{Z}_{x_0}\|_\alpha+\|Y\|_{B^\beta_{4,\infty}}\|:\bar{Z}_{x_0}^2:\|_\alpha]ds
\\\leq& C\int_0^1(1-s)^{-\frac{\beta+1/2-\alpha}{2}}[\|Y^2\|_{B^{\beta+1/2}_{2,\infty}}\|\bar{Z}_{x_0}\|_\alpha+\|Y\|_{B^{\beta+1/2}_{2,\infty}}\|:\bar{Z}_{x_0}^2:\|_\alpha]ds
\\\leq&C(\|\underline{Z}\|_{\mathfrak{L}_1}, \|x_0\|_{\alpha})\int_0^1(1-s)^{-\frac{\beta+1/2-\alpha}{2}}[(\|\nabla Y^2\|_{L^2}^{\beta+1/2+\varepsilon}\|Y^2\|^{1/2-\beta-\varepsilon}_{L^2}+\|Y^2\|_{L^2})
\\&+s^{-\beta-\alpha}(\|\nabla Y\|_{L^2}^{\beta+1/2+\varepsilon}\|Y\|^{1/2-\beta-\varepsilon}_{L^2}+\|\nabla Y\|_{L^2})]ds
\\\leq& C(\|\underline{Z}\|_{\mathfrak{L}_1}, \|x_0\|_{\alpha})+C(\|\underline{Z}\|_{\mathfrak{L}_1}, \|x_0\|_{\alpha})
\int_0^1(\|\nabla Y^2\|_{L^2}^{2}+\|Y^2\|_{L^2})
ds\\\leq& C(\|\underline{Z}\|_{\mathfrak{L}_1}, \|x_0\|_{\alpha}),\endaligned$$
where $0<\varepsilon<\frac{1}{2}-\beta$ and  we used Lemma 2.3 in the first inequality,  we used Lemma 2.1 in the second inequality, and the fact that
$$\|\cdot\|_{B^{\beta+\frac{1}{2}}_{2,\infty}}\leq C\|\cdot\|_{B^{\beta+\frac{1}{2}}_{2,1}}\leq C\|\cdot\|_{H^{\beta+\frac{1}{2}+\varepsilon}_{2}},$$
and Lemma 2.4 in the third inequality, and we used (3.7) for $p=2$ and $4$ in the last two inequalities. Combining the above estimates (3.8) follows.

Combining (3.5) and (3.8) and using $\|Y(1)\|_{L^p}\leq C\|Y(1)\|_{\beta}$, the  result follows.
 $\hfill\Box$
\vskip.10in

The proof of Lemma 3.5 is similar to that of Lemma 3.4. But we should pay attention to how each term depends on $\lambda$, as $\tilde{Y}$ depends on $\lambda$.

\vskip.10in

\th{Lemma 3.5} (i) For every even $p>1$, there exist constants $C(p), C(p,\lambda,\|\underline{Z}\|_{\mathfrak{L}_1},\|x_0\|_\alpha,\|x_1\|_\alpha)>0, \gamma(p)>1$ independent of $\omega, t$ such that for every $t\geq1$ and $\omega\in\Omega_0$
$$\int_1^t\|\tilde{Y}(s)\|_{L^p}^pds\leq  C(p)\int_{1}^t (1+\sum_{i=0}^1\|x_i\|_\alpha^{\gamma(p)}+\|Z\|_{\alpha}^{\gamma(p)}+\sum_{n=2}^3\|:Z^n:\|_{\alpha}^{{\gamma(p)}})ds+C(p,\lambda,\|\underline{Z}\|_{\mathfrak{L}_1},\|x_0\|_\alpha,\|x_1\|_\alpha).$$
$$\aligned&\int_1^t\|\nabla \tilde{Y}(s)\|_{L^2}^2ds+\|\tilde{Y}(t)\|_{L^p}^p\\\leq&  C(p)\lambda\int_{1}^t (1+\sum_{i=0}^1\|x_i\|_\alpha^{\gamma(p)}
+\|Z\|_{\alpha}^{\gamma(p)}+\sum_{n=2}^3\|:Z^n:\|_{\alpha}^{\gamma(p)})ds+C(p,\lambda,\|\underline{Z}\|_{\mathfrak{L}_1},\|x_0\|_\alpha,\|x_1\|_\alpha).\endaligned$$

\proof
\textbf{Step 1} We first prove that for every even $p>1$ there exists $\gamma(p)>1$ such that
\begin{equation}\aligned\int_{1}^t\|\tilde{Y}(s)\|_{L^p}^pds\leq& \|\tilde{Y}(1)\|_{L^p}^p+ C(p)\int_{1}^t (1+\sum_{i=0}^1\|x_i\|_\alpha^{\gamma(p)}
+\|Z\|_{\alpha}^{\gamma(p)}+\sum_{n=2}^3\|:Z^n:\|_{\alpha}^{\gamma(p)})ds\\&+C(\|\underline{Z}\|_{\mathfrak{L}_1}, \|x_0\|_{\alpha})\endaligned.\end{equation}

Similarly as in the proof of Lemma 3.4
we have that for $t\geq1$ and even $p>1$
\begin{equation}\aligned&\frac{1}{p}\|\tilde{Y}(t)\|_{L^p}^p+\lambda\int_1^t\|\tilde{Y}(s)\|_{L^p}^pds
+\int_1^t[(p-1)\langle\nabla \tilde{Y}(s), \tilde{Y}(s)^{p-2}\nabla \tilde{Y}(s)\rangle+a_1\|\tilde{Y}(s)^{p+2}\|_{L^1}]ds\\
=&-\int_1^t[\|\tilde{Y}(s)\|_{L^p}^p+\langle \Psi(\tilde{Y}(s),\bar{Z}_{x_1}(s)),\tilde{Y}(s)^{p-1}\rangle]ds
+\lambda\int_1^t\langle Y(s),\tilde{Y}(s)^{p-1}\rangle ds\\&-\lambda\int_1^t\langle e^{sA}(x_1-x_0),\tilde{Y}(s)^{p-1}\rangle ds
+\frac{1}{p}\|\tilde{Y}(1)\|_{L^p}^p.\endaligned\end{equation}

Now by similar calculations as in the proof of Lemma 3.4 and using (3.7) we  deduce that there exist $ \gamma(p)>1$ such that
\begin{equation}\aligned&\|\tilde{Y}(t)\|_{L^p}^p+\lambda\int_1^t\|\tilde{Y}(s)\|_{L^p}^pds+\frac{p-1}{2}\int_1^t\||\nabla\tilde{Y}(s)|^2\tilde{Y}^{p-2}(s)\|_{L^1}ds
\\\leq&C(p)\int_1^t (1+\|\bar{Z}_{x_1}\|_{\alpha}^{\gamma_0}+\sum_{n=2}^3\|:\bar{Z}^n_{x_1}:\|_{\alpha}^{p+2})ds+\lambda p\int_1^t\|Y(s)\|_{L^p}
\|\tilde{Y}(s)\|_{L^p}^{p-1}ds\\&+C\lambda p\int_1^t s^{-\frac{\beta-\alpha}{2}}\|x_0-x_1\|_{\alpha}
\|\tilde{Y}(s)\|_{L^p}^{p-1}ds+\|\tilde{Y}(1)\|_{L^p}^p\\\leq&C(p)\int_1^t (1+\|{Z}\|_{\alpha}^{\gamma(p)}+\sum_{n=2}^3\|:{Z}^n:\|_{\alpha}^{\gamma(p)}
+\|x_1\|_\alpha^{\gamma(p)})ds+\lambda C(p)\int_1^t\|Y(s)\|_{L^p}^pds+
\frac{\lambda}{2}\int_1^t\|\tilde{Y}(s)\|_{L^p}^pds\\&+C(p)\lambda\int_1^t \|x_0-x_1\|_{\alpha}^p
ds+\|\tilde{Y}(1)\|_{L^p}^p,\endaligned\end{equation}
where  we used H\"{o}lder's inequality and Lemma 2.2 to control $\|e^{sA}(x_0-x_1)\|_{\mathcal{C}^\beta}$ by $Cs^{-\frac{\beta-\alpha}{2}}\|x_0-x_1\|_{\mathcal{C}^\alpha}$ for $\beta>-\alpha$
in the first inequality and we used Young's inequality in the second inequality.
By  Lemma 3.4 and the fact that $\lambda>1$, (3.9) follows.

\textbf{Step 2} We prove that for even $p>1$ with $(-2\alpha+\delta)(p+2)<1$
\begin{equation}\sup_{0\leq t\leq 1}\|\tilde{Y}(t)\|_{L^p}^p+\int_0^1\|\tilde{Y}^{p-2}|\nabla \tilde{Y}|^2\|_{L^1}ds\leq \lambda C(p, \|\underline{Z}\|_{\mathfrak{L}_1}, \|x_0\|_{\alpha}, \|x_1\|_{\alpha}).\end{equation}
By similar arguments as in Step 1 we have for $0\leq t\leq 1$, even $p>1$ with $-2\alpha(p+2)<1$ and $\varepsilon>0$ small enough
$$\aligned&\|\tilde{Y}(t)\|_{L^p}^p+\lambda\int_0^t\|\tilde{Y}(s)\|_{L^p}^pds+\frac{p-1}{2}\int_0^t\||\nabla\tilde{Y}(s)|^2\tilde{Y}^{p-2}(s)\|_{L^1}ds
\\\leq&C(p)\int_0^t (1+\|\bar{Z}_{x_1}\|_{\alpha}^{\gamma_0}+\sum_{n=2}^3\|:\bar{Z}^n_{x_1}:\|_{\alpha}^{p+2})ds+\lambda\int_0^t\|Y(s)\|_{L^p}
\|\tilde{Y}(s)\|_{L^p}^{p-1}ds\\&+C\lambda p\int_0^t s^{2\alpha}\|x_0-x_1\|_{\alpha}
\|\tilde{Y}(s)\|_{L^p}^{p-1}ds\\\leq&C(p,\|\underline{Z}\|_{\mathfrak{L}_1})\int_0^t s^{(2\alpha-\varepsilon)(p+2)}(1+\|x_1\|_\alpha^{\gamma(p)})ds
+\lambda C(p)\int_0^t\|Y(s)\|_{L^p}^pds+
\frac{\lambda}{2}\int_0^t\|\tilde{Y}(s)\|_{L^p}^pds\\&+C(p)\lambda\int_0^ts^{2\alpha p} \|x_0-x_1\|_{\alpha}^p
ds,\endaligned$$
where we used Lemma 2.7 in the last inequality. By  Lemma 3.4, (3.7) and the fact that $\lambda>1$, (3.12) follows.

\textbf{Step 3} We prove that for $0<\beta<\frac{1}{2}+\alpha$
\begin{equation}\|\tilde{Y}(1)\|_{\beta}\leq C(\lambda,\|\underline{Z}\|_{\mathfrak{L}_1}, \|x_0\|_{\alpha}, \|x_1\|_{\alpha}).\end{equation}
Since $\tilde{Y}$ satisfies the mild equation, by similar arguments as in Step 3 in the proof of Lemma 3.4 we have
$$\aligned\|\tilde{Y}(1)\|_{\beta}\leq&C\int_0^1(1-s)^{-\frac{\beta}{2}-\frac{1}{2}}[\|\tilde{Y}^3\|_{L^2}+\lambda(\|\tilde{Y}\|_{L^2}+\|Y\|_{L^2})]ds
+C\int_0^1(1-s)^{-\frac{\beta+1/2-\alpha}{2}}[\|\tilde{Y}^2\bar{Z}_{x_1}\|_{B^\alpha_{4,\infty}}\\&+\|\tilde{Y}:\bar{Z}_{x_1}^2:\|_{B^\alpha_{4,\infty}}]ds
+C\int_0^1(1-s)^{-(\beta-\alpha)/2}(\|:\bar{Z}_{x_1}^3:\|_\alpha+\lambda\|x_0\|_\alpha+\lambda\|x_1\|_\alpha) ds
\\\leq& C(p,\lambda, \|\underline{Z}\|_{\mathfrak{L}_1}, \|x_0\|_{\alpha}),\endaligned$$
where for the term in the first integral we used (3.12) and for the term in the second and third integral we used similar arguments as in Step 3 in the proof of Lemma 3.4.
Combining (3.9) and (3.13) the first result follows.

The second  follows from (3.11), (3.12) and Lemma 3.4.
 $\hfill\Box$
\vskip.10in

In the following we give an estimate of the Wick power $:Z^k:$, which is required in the proof of the main results.
\vskip.10in
For  $\gamma>0$ and $K>0$ we introduce the following notations:
\begin{equation}E_{K,\gamma}:=\{\|\bar{Z}\|_{\mathfrak{L}_1}\leq K, \int_1^t[\|Z\|_{\alpha}^{\gamma}
+\|:Z^2:\|_{\alpha}^{\gamma}+\|:Z^3:\|_{\alpha}^\gamma ]ds\leq K(1+t), \forall t\geq1 \}.\end{equation}

\vskip.10in
\th{Lemma 3.6} For every $\gamma>0, \varepsilon>0$ there exists $K>0$ such that $P(E_{K,\gamma})\geq 1-\varepsilon$.

\proof To prove this result, we first introduce the following stationary Markov process. Define $Z_1(t)=\int_{-\infty}^te^{(t-s)A}dW(s)$. We  also define
$$:Z_1^2::=\lim_{\varepsilon\rightarrow0}[(Z_1*\rho_\epsilon)^2-c_\epsilon]\textrm{ in } L^p(\Omega,C([0,\infty),
\mathcal{C}^\alpha)),$$ $$:Z_1^3::=\lim_{\epsilon\rightarrow0}[(Z_1*\rho_\epsilon)^3-3c_\epsilon Z_1*\rho_\epsilon]\textrm{ in }
L^p(\Omega,C([0,\infty),\mathcal{C}^\alpha)),$$ for $p>1$ as in [RZZ15, Lemma 3.3], which are also stationary Markov processes. Here $\rho_\epsilon$ and $c_\epsilon$ are introduced in Section 2.2.
By [DZ96, Theorem 3.3.1] we know that for every $q>1$ there exists $\eta\in L^2(\Omega,P)$ such that
$$\mathfrak{Z}_T:=\frac{1}{T}\int_0^T[\|Z_1\|_{\alpha}^q+\sum_{n=2}^3\|:Z^n_1:\|_{\alpha}^q] ds\rightarrow \eta, \textrm{ as }T\rightarrow\infty,\quad P-a.s.,$$
which implies that for every $\varepsilon>0$, there exists $\Omega_1\subset \Omega$ such that $P(\Omega_1)<\varepsilon/4$
and $$\sup_{\omega\in \Omega_1^c}|\mathfrak{Z}_T(\omega)-\eta(\omega)|\rightarrow0, \textrm{ as }T\rightarrow\infty.$$
From this we can deduce that there exists $ T_0$ independent of $\omega$ such that for $T\geq T_0$
$$\mathfrak{Z}_T(\omega)\leq \eta(\omega)+1, \forall \omega\in \Omega_1^c,$$
which combined with $\eta\in L^2(\Omega;P)$ yields that there exists $K_1>0$ such that
$$P\{\int_0^T[\|Z_1\|_{\alpha}^{2\gamma}+\sum_{n=2}^3\|:Z^n_1:\|_{\alpha}^{2\gamma}] ds\leq K_1T, \forall T\geq T_0\}>1-\varepsilon/3.$$
Thus, there
 exists $K_2>0$ such that
$$P\{\int_0^T[\|Z_1\|_{\alpha}^{2\gamma}+\sum_{n=2}^3\|:Z^n_1:\|_{\alpha}^{2\gamma}] ds\leq K_2(T+1), \forall T\geq0\}>1-\varepsilon/3.$$
Now we give the relations of $Z$ and $Z_1$.
By (2.2) and similar arguments as in the proof of [RZZ15, Lemma 3.6] we have that
$$\aligned Z(t)=&Z_1(t)-e^{tA}Z_1(0),\\:Z^2(t):=&:Z_1^2(t):-2e^{tA}Z_1(0)Z_1(t)+(e^{tA}Z_1(0))^2,
\\:{Z}^3(t):=&:Z_1^3(t):+3(e^{tA}Z_1(0))^2Z_1(t)-3e^{tA}Z_1(0):Z_1(t)^2:-(e^{tA}Z_1(0))^3,\endaligned$$
which combined with Lemma 2.3 implies that for $\beta>-\alpha>0$
$$\aligned\|Z(t)\|_{\alpha}\leq& \|Z_1(t)\|_{\alpha}+\|Z_1(0)\|_{\alpha},\\\|:Z^2(t):\|_{\alpha}\leq& C[ \|:Z_1^2(t):\|_{\alpha}+\|e^{tA}Z_1(0)\|_{\beta}\|Z_1(t)\|_{\alpha}+\|e^{tA}Z_1(0)\|_{\beta}^2],
\\\|:Z^3(t):\|_{\alpha}\leq& C[ \|:Z_1^3(t):\|_{\alpha}+\|e^{tA}Z_1(0)\|_{\beta}^2\|Z_1(t)\|_{\alpha}
\\&+\|e^{tA}Z_1(0)\|_{\beta}\|:Z_1^2(t):\|_{\alpha}+\|e^{tA}Z_1(0)\|_{\beta}^3].\endaligned$$
Now using Lemma 2.2 we have
$$\aligned&\int_1^T[\|Z\|_{\alpha}^{\gamma}+\|:Z^2:\|_{\alpha}^{\gamma}+\|:Z^3:\|_{\alpha}^{\gamma}] ds\\\leq& C\int_1^T[\|Z_1(0)\|_{\alpha}^{4\gamma}+1+\|Z_1(s)\|_{\alpha}^{2\gamma}
+\|:Z^2_1(s):\|_{\alpha}^{2\gamma}+\|:Z^3_1(s):\|_{\alpha}^\gamma] ds,\endaligned$$
which implies that there exists $K_3>0$ such that
$$P\{\int_1^T[\|Z\|_{\alpha}^{\gamma}+\sum_{n=2}^3\|:Z^n:\|_{\alpha}^{\gamma}] ds\leq K_3(T+1), \forall T\geq0\}>1-\varepsilon/2.$$
On the other hand, by Lemma 2.6 we have $E\|\underline{Z}\|_{\mathfrak{L}^1}^2<\infty$, which implies that there exist $K_4>0$ such that
$$P(\|\underline{Z}\|_{\mathfrak{L}^1}\leq K_4)>1-\varepsilon/2.$$
Combining the above results we obtain that
there exists $K>0$ such that
$$P(E_{K,\gamma})\geq 1-\varepsilon.$$
 $\hfill\Box$

\section{Uniqueness of the invariant measure}

In this section we will prove our main results: uniqueness of the invariant measure.
We first present an abstract result based on asymptotic coupling from [HMS11]: Let $\mathcal{P}$ be a Markov transition function on a Polish space $(\mathbb{X},\rho)$ and let $\mathbb{X}_\infty= \mathbb{X}^{\mathbb{N}}$ be the associated space of one-sided infinite
sequences. Denote the collection of all Borel probability measures on $\mathbb{X}$ by $\mathcal{M}(\mathbb{X})$. Take
$\mathcal{P}_\infty:\mathbb{X}\rightarrow\mathcal{M}(\mathbb{X}_\infty)$ to be the probability kernel defined by stepping with the Markov kernel
$\mathcal{P}$. For $\mu_0\in\mathcal{M}(\mathbb{X})$, let $\mu_0 \mathcal{P}_\infty\in \mathcal{M}(\mathbb{X}_\infty)$ be the measure defined by
$\int_{\mathbb{X}}\mathcal{P}_\infty(x,\cdot)d\mu_0(x)$. Given $\mu_1, \mu_2\in\mathcal{M}(\mathbb{X})$, consider
$$\tilde{\mathcal{C}}(\mu_1\mathcal{P}_\infty,\mu_2\mathcal{P}_\infty):=\{\Gamma\in \mathcal{M}(\mathbb{X}_\infty\times\mathbb{X}_\infty):
\Pi_i^\sharp\Gamma\ll \mu_i\mathcal{P}_\infty\textrm{ for each } i\in \{1,2\}\}, $$
where $\Pi_i$ is the projection onto the ith coordinate and $f^\sharp\mu_0$ is the push-forward of the measure $\mu_0$ defined by $(f^\sharp\mu_0)(B)=\mu_0(f^{-1}(B))$.
We also denote the diagonal at infinity by
$$D:=\{(x,y)\in \mathbb{X}_\infty\times \mathbb{X}_\infty:\lim_{n\rightarrow\infty}\rho(x_n,y_n)=0\}.$$
\vskip.10in

\th{Theorem 4.1} ([HMS11, Corollary 2.2]). Suppose that there exists a Borel measurable set $B\subset \mathbb{X}$ such that

(i) for any $\mathcal{P}$ invariant Borel probability measure $\mu$, $\mu(B)>0$,

(ii) there exists a measurable map $\Gamma: B\times B\rightarrow \mathcal{M}(\mathbb{X}_\infty\times \mathbb{X}_\infty)$ such that, for all $x,y\in B$, $\Gamma_{x,y}\in \tilde{\mathcal{C}}(\delta_x\mathcal{P}_\infty,\delta_y\mathcal{P}_\infty)$ and $\Gamma_{x,y}(D)>0$,

then there exists at most one invariant probability measure for $\mathcal{P}$.

\vskip.10in

Now we prove our main result by using Theorem 4.1.

\vskip.10in
\no \emph{Proof of Theorem 1.1.}  For  a given stochastic basis $(\Omega,\mathcal{F},(\mathcal{F}_t)_{t\geq0},P)$ and a cylindrical Wiener process $W$ as in Section 3, we use $X(x_0)$ and $\tilde{X}(x_1)$ to denote the solutions of the equations (3.2) and (3.4) obtained in Section 3 starting from $x_0, x_1\in\mathcal{C}^\alpha$, respectively.  We choose $B=\mathbb{X}=\mathcal{C}^\alpha$. Let $P_\infty:\mathcal{C}^\alpha\mapsto\mathcal{M}(\mathcal{C}^\alpha_\infty)$ be the probability kernel corresponding to $X$ evaluated at integer times.

\textbf{Girsanov transform }

Set $v=\lambda(\tilde{X}(x_1)-X(x_0))$ and let $\tilde{W}(t)=W(t)-\int_0^{t\wedge\tau_R}v(s)ds$, where
$$\tau_R:=\inf\{t>0, \int_0^t\|X(x_0,s)-\tilde{X}(x_1,s)\|_{L^2}^2ds\geq R\}.$$
Since
$$E\exp\bigg(\frac{1}{2}\int_0^{\tau_R}\|v(s)\|_{L^2}^2ds\bigg)\leq e^{\frac{1}{2}R\lambda^2},$$
by the Girsanov theorem there is a probability measure $Q$ on $(\Omega,\mathcal{F},(\mathcal{F}_t)_{t\geq0})$ such that under $Q$, $\tilde{W}$ is a standard Wiener process. Moreover, it holds that $P\thicksim Q$ on $\mathcal{F}_\infty=\sigma(\cup_{t\geq0}\mathcal{F}_t)$.

\textbf{Construction of the coupling }

Let $\hat{Z}$ is the solution to the following linear equation
$$d\hat{Z}(t)=A\hat{Z}(t)dt+d\tilde{W}(t),\quad \hat{Z}(0)=0,$$
and $:\hat{Z}^k:$ can be defined similarly as $:Z^k:$ as in Section 2. Moreover, we use similar notations as in Section 2: $\bar{\hat{Z}}_{x_1}:=\hat{Z}+e^{tA}x_1$, and for $n=2, 3$,
$$:\bar{\hat{Z}}_{x_1}^n(t)::=\sum_{k=0}^nC_n^k(e^{tA}x_1)^{n-k}:\hat{Z}^k(t):.$$

Furthermore, we derive the relation between different Wick powers under $P$ and $Q$ respectively.  Since $\hat{Z}=Z+a$ with $a(t)=-\int_0^te^{(t-s)A}1_{s\leq \tau_R}v(s)ds\in C([0,\infty);\mathcal{C}^\beta)$ for some $\beta>-\alpha$, we have that there exists $\Omega_0'\subset\Omega_0$ such that $P(\Omega_0')=1$ and the following hold for $\omega\in\Omega_0'$ in $C((0,\infty);\mathcal{C}^\alpha)$
\begin{equation}:\hat{Z}^2:=:Z^2:+2Za+a^2,\end{equation}
and
\begin{equation}:\hat{Z}^3:=:Z^3:+3Za^2+3:Z^2:a+a^3.\end{equation}
We will prove (4.1) and (4.2) at the end of the proof.

For $\omega\in\Omega_0'$, by (4.1), (4.2) and $a\in  C([0,\infty);\mathcal{C}^\beta)$ we know that for $n=2,3, T\in\mathbb{R}^+$
$$\hat{Z}\in C([0,T];{\mathcal{C}}^{\alpha}), :\hat{Z}^n:\in C((0,T]; {\mathcal{C}}^{\alpha}) ,\|\underline{{\hat{Z}}}\|_{\mathfrak{L}_T}<\infty,$$
 which by Theorem 3.1 implies that for $\omega\in\Omega_0'$ there exists a unique mild solution $\hat{Y}(\omega)\in C([0,\infty),\mathcal{C}^\beta)$
 to the following equation
\begin{equation}\frac{d\hat{Y}}{dt}=A\hat{Y}-a_1\hat{Y}^3+\Psi(\hat{Y},{\bar{\hat{Z}}_{x_1}}),\quad \hat{Y}(0)=0.\end{equation}
Define  $$\hat{X}(x_1,\omega)=\left\{\begin{array}{ll}[\hat{Y}+e^{tA}x_1+\hat{Z}](\omega),&\ \ \ \ \textrm{ if } \omega\in\Omega_0'\\0&\ \ \ \ \textrm{ otherwise } .\end{array}\right.$$
Then we conclude that under $Q$, $\hat{X}(x_1)$ is also a mild solution to the equation (3.2) with $x=x_1$ and with $W$ replaced by $\tilde{W}$, which combined with Theorem 3.1 and the Yamada-Watanabe Theorem in [Kurz07] implies that under $Q$, $\hat{X}(x_1)$ has the same law as the solution $X(x_1)$ to the equation (3.2) starting from $x_1$. Since $P\thicksim Q$, we have
that under $P$ the law of the pair $(X(x_0),\hat{X}(x_1))$ has marginals
which are equivalent  to the marginals of the solutions to (3.2) starting respectively from $x_0$ and $x_1$.  Set $\Gamma_{x_0,x_1}:=$ law of $(X(x_0),\hat{X}(x_1))$ for $x_0,x_1\in\mathcal{C}^\alpha$. It follows that $\Gamma_{x_0,x_1}\in \tilde{\mathcal{C}}(\delta_{x_0}P_\infty,\delta_{x_1}P_\infty)$. It remains to show that $\Gamma_{x_0,x_1}(D)>0$.

We have that $\hat{X}(x_1)$ satisfies the following equation in the mild sense $P$-a.s.:
$$\aligned d\hat{X}=&[A\hat{X}-a_1(\hat{X}-\bar{\hat{Z}}_{x_1})^3+\Psi(\hat{X}-\bar{\hat{Z}}_{x_1},{\bar{\hat{Z}}_{x_1}})]dt+d\tilde{W}.\endaligned$$
By (4.1), (4.2) we have that there exists $\Omega_2\subset \Omega_0'$ such that $P(\Omega_2)=1$ and for $\omega\in\Omega_2$,
 $\hat{X}(x_1,\omega)$ also satisfies the following equation in the mild sense :
$$d\hat{X}=[A\hat{X}-a_1(\hat{X}-\bar{{Z}}_{x_1})^3+\Psi(\hat{X}-\bar{{Z}}_{x_1},{\bar{{Z}}_{x_1}})]dt+dW-v1_{t\leq \tau_R}dt.$$
Then on $\{\tau_R=\infty\}\cap \Omega_2$, $\hat{X}-\bar{Z}_{x_1}$ also satisfies (3.3). By Theorem 3.3 we obtain that $\hat{X}-\bar{Z}_{x_1}=\tilde{Y}$ on $\{\tau_R=\infty\}\cap\Omega_2$, which implies that
$\hat{X}=\tilde{X}$ on $\{\tau_R=\infty\}\cap\Omega_2$.
Here $\tilde{Y}$ is the solution to (3.3) and $\tilde{X}(x_1)=\tilde{Y}+e^{tA}x_1+Z$. Now to prove $\Gamma_{x_0,x_1}(D)>0$, it suffices to estimate $X(x_0)-\tilde{X}(x_1)$.

\textbf{Estimate of $X(x_0)-\tilde{X}(x_1)$}:

In the following we  estimate $X(x_0)-\tilde{X}(x_1)$ and we do all the calculations informally, but all the calculations below can
be made rigorous by approximation as done in the proof of [RZZ15, Theorem 3.10].  Set  $Y_1=Y+e^{tA}x_0$, $\tilde{Y}_1=\tilde{Y}+e^{tA}x_1$ and
$u=X-\tilde{X}=Y_1-\tilde{Y}_1$.
By the binomial formula (2.2) we have that $P$-a.s. $Y_1$ and $\tilde{Y}_1$ are the mild solutions to the following equations
$$\frac{d}{dt}Y_1=AY_1-[a_1Y_1^3+\Psi(Y_1,{Z})],\quad Y_1(0)=x_0,$$
and
$$\frac{d}{dt}\tilde{Y}_1=A\tilde{Y}_1+\lambda (Y_1-\tilde{Y}_1)-[a_1\tilde{Y}_1^3+\Psi(\tilde{Y}_1,{Z})],\quad Y_1(0)=x_1$$
respectively.
It is obvious  that $P$-a.s. $u$ is the mild solution to the following equation:
$$\frac{d}{dt}u=Au-\lambda u-[a_1Y_1^3-a_1\tilde{Y}_1^3+\Psi(Y_1,{Z})-\Psi(\tilde{Y}_1,{Z})],\quad u(0)=x_0-x_1.$$

Standard energy estimates yield
\begin{equation}\frac{1}{2}\frac{d}{dt}\|u\|_{L^2}^2+\|\nabla u\|_{L^2}^2+\lambda \|u\|_{L^2}^2\leq-\langle(\Psi(Y_1,\underline{Z})-\Psi(\tilde{Y}_1,{Z})),u\rangle,\end{equation}
where we used that $$-\langle Y_1^3-\tilde{Y}_1^3,u\rangle\leq 0.$$
Now we calculate each term in $\langle(\Psi(Y_1,Z)-\Psi(\tilde{Y}_1,Z)),u\rangle$: For the first term we have
\begin{equation}\aligned&3a_1|\langle (Y_1^2-\tilde{Y}_1^2)Z,u\rangle|=3a_1|\langle u^2,(Y_1+\tilde{Y}_1)Z\rangle|\\\leq&C_S \|u^2\|_{B^{-\alpha}_{\frac{4}{3},1}}(\|Y_1Z\|_{B^{\alpha}_{4,\infty}}+\|\tilde{Y}_1Z\|_{B^{\alpha}_{4,\infty}})
\\\leq& C_S\|\Lambda^{\frac{1}{2}}u\|_{L^2}^2(\|Y_1Z\|_{B^{\alpha}_{4,\infty}}+\|\tilde{Y}_1Z\|_{B^{\alpha}_{4,\infty}})\\\leq &C_S\|u\|_{L^{2}}
(\|\nabla u\|_{L^{2}}+\|u\|_{L^{2}})(\|Y_1Z\|_{B^{\alpha}_{4,\infty}}+\|\tilde{Y}_1Z\|_{B^{\alpha}_{4,\infty}})\\\leq &C_S\|u\|_{L^{2}}^{2}
(\|Y_1Z\|_{B^{\alpha}_{4,\infty}}^2+\|\tilde{Y}_1Z\|_{B^{\alpha}_{4,\infty}}^2+1)+\frac{1}{4}\|\nabla u\|_{L^{2}}^2,\endaligned\end{equation}
where $C_S$ is a constant changing from line to line and we used Lemma 2.3 in the first inequality and  Lemmas 2.1 and 2.4 to deduce that
\begin{equation}\|u^2\|_{B^{-\alpha}_{\frac{4}{3},1}}\leq C_S\|\Lambda^s u^2\|_{L^{\frac{4}{3}}}\leq C_S\|\Lambda^s u\|_{L^{2}}\|u\|_{L^{4}}\leq C_S\|\Lambda^{\frac{1}{2}}u\|_{L^2}^2,\end{equation}
for $\frac{1}{2}>s>-\alpha$ in the second inequality and we used Lemma 2.4 in the third inequality and Young's inequality in the last inequality.
For the second term we have
\begin{equation}\aligned 3a_1|\langle Y_1:Z^2:-\tilde{Y}_1:Z^2:, u\rangle|
\leq& C_S\|u^2\|_{{B^{-\alpha}_{1,1}}}\|:Z^2:\|_{\alpha}\\\leq &C_S\|u\|_{L^{2}}
(\|\nabla u\|_{L^{2}}+\|u\|_{L^2})\|:Z^2:\|_{\alpha}\\\leq &C_S\|u\|_{L^{2}}^{2}(\|:Z^2:\|_{\alpha}^2+1)
+\frac{1}{4}\|\nabla u\|_{L^{2}}^2,\endaligned\end{equation}
where we used Lemma 2.4 in the first inequality,
 Lemmas 2.1 and (4.6) to deduce that
$$\|u^2\|_{B^{-\alpha}_{1,1}}\leq\|u^2\|_{B^{-\alpha}_{\frac{4}{3},1}}\leq C_S\|\Lambda^{\frac{1}{2}}u\|_{L^2}^2,$$
for $\frac{1}{2}>s>-\alpha, q>1, \frac{1}{q}=\frac{1}{q_1}+\frac{1}{q_2}, q_2<4, \frac{2}{q_1}-s>\frac{1}{2}$ in the second inequality
 and we used Young's inequality in the last inequality.

For the last term
we have
\begin{equation}\aligned |a_2\langle Y_1-\tilde{Y}_1, u\rangle|
\leq C\|u\|_{L^{2}}^{2}.\endaligned\end{equation}

Combining (4.4)-(4.8) we obtain
$$\aligned &\frac{1}{2}\frac{d}{dt}\|u\|_{L^2}^2+\lambda \|u\|_{L^2}^2\\\leq& \|u\|_{L^{2}}^{2}
C_S[\|Y_1{Z}\|_{B^\alpha_{4,\infty}}^2+\|\tilde{Y}_1{Z}\|_{B^\alpha_{4,\infty}}^2+\|:Z^2:\|_{\alpha}^2+1]:=\|u\|_{L^{2}}^{2}L.\endaligned$$

Then Gronwall's inequality yields that
$$\|u(t)\|_{L^2}^2\leq \|u(1)\|_{L^2}^2\exp {\int_1^t 2(-\lambda+L(s))ds}.$$
Here we use Gronwall's inequality starting from $t=1$ instead of $t=0$ since $u(0)$ is not in $L^2$.

Recall that for  $\gamma, K>0$, $E_{K,\gamma}$ has been defined in (3.14).
By Lemma 3.6 we know that for every $\gamma>0$  there exists  $K>0$, such that  $P(E_{K,\gamma})>0$.
In the following we estimate each term in $L$ on $E_{K,\gamma}$ with  $\gamma>0$ to be determined later: We have that on $E_{K,\gamma}$
\begin{equation}\aligned\int_1^t\|Y_1Z\|_{B^\alpha_{4,\infty}}^{2}ds\leq& C_S\int_1^t\| Y_1\|_{B^{\beta}_{4,\infty}}^{2}\|Z\|_{\alpha}^{2}ds
\\\leq& C_S\int_1^t(\|\nabla Y_1\|_{L^2}^{2\beta_0}\| Y_1\|_{L^2}^{2(1-\beta_0)}+\| Y_1\|_{L^2}^2)\|Z\|_{\alpha}^{2}ds\\\leq& C_S[(\int_1^t\|\nabla Y_1\|_{L^2}^{2\beta_0 p_1}ds)^{\frac{1}{p_1}} (\int_1^t\| Y_1\|_{L^2}^{2(1-\beta_0)p_2}ds)^{\frac{1}{p_2}}(\int_1^t\|Z\|_{\alpha}^{2p_3}ds)^{\frac{1}{p_3}}\\&+(\int_1^t\| Y_1\|_{L^2}^{4}ds)^{\frac{1}{2}}(\int_1^t\|Z\|_{\alpha}^{4}ds)^{\frac{1}{2}}],\endaligned\end{equation}
where $\frac{1}{p_1}+\frac{1}{p_2}+\frac{1}{p_3}=1, p_i>1, i=1,2,3$, $\beta>-\alpha$, $\beta_0=\beta+\frac{1}{2}+\varepsilon, \varepsilon>0, 2\beta_0 p_1\leq 2$,  and we used Lemma 2.3 in the first inequality,  Lemma 2.1 to deduce that
$$\|Y_1\|_{B^\beta_{4,\infty}} \leq C_S\|Y_1\|_{B^{\beta+1/2}_{2,\infty}}\leq C_S\|Y_1\|_{B^{\beta+1/2}_{2,1}}\leq C_S\|Y_1\|_{H^{\beta_0}_{2}}$$ in the second inequality and  H\"{o}lder's inequality in the last inequality. In the following we estimate each term on the right hand side of (4.9):
By  Lemma 3.4 we know that for any even $p>1$ we have on $E_{K,\gamma}$ with $\gamma\geq\gamma(p)$
$$\aligned\int_1^t\|Y_1(s)\|_{L^p}^pds\leq& C(p)[\int_1^t\|Y(s)\|_{L^p}^pds+\int_1^t\|e^{sA}x_0\|_{L^p}^pds]\\\leq& C(p)\int_1^t [1+\|Z\|_{\alpha}^{\gamma(p)}+\|:Z^2:\|_{\alpha}^{\gamma(p)}+\|:Z^3:\|_{\alpha}^{\gamma(p)}]ds+C(\|x_0\|_{\alpha})(1+t)
\\&+C(\|\underline{Z}\|_{\mathfrak{L}_1}, \|x_0\|_{\alpha})\\\leq &C(p, \|x_0\|_{\alpha},K)(1+t),\endaligned$$
where we used Lemma 2.2 to control $\|e^{sA}x_0\|_{L^p}\leq C_Ss^{\alpha }\|x_0\|_\alpha$ in the second inequality.
Similarly, by Lemma 3.4 we have that on $E_{K,\gamma}$ for $\gamma\geq\gamma(2)$
$$\aligned\int_1^t\|\nabla Y_1(s)\|_{L^2}^{2\beta_0 p_1}ds\leq& C(p_1)[\int_1^t\|\nabla Y(s)\|_{L^2}^2ds+t+\int_1^t\|\nabla e^{sA}x_0\|_{L^2}^{2\beta_0 p_1}ds]\\\leq& C(p_1)\int_1^t [1+\|x_0\|_\alpha^{\gamma(2)}+\|Z\|_{\alpha}^{\gamma(2)}+\|:Z^2:\|_{\alpha}^{\gamma(2)}+\|:Z^3:\|_{\alpha}^{\gamma(2)}]ds+t\\&+C\int_1^ts^{-(1+\varepsilon-\alpha)
\beta_0 p_1}\|x_0\|_{\alpha}^{2\beta_0 p_1}ds+C(\|\underline{Z}\|_{\mathfrak{L}_1}, \|x_0\|_{\alpha})\\\leq &C(p_1, \|x_0\|_{\alpha},K)(1+t),\endaligned$$
where  we used Young's inequality and $2\beta_0p_1\leq 2$ in the first inequality and Lemmas 2.1, 2.2 to deduce that $\|\nabla e^{sA}x_0\|_{L^2}\leq C_Ss^{-(1+\varepsilon-\alpha)/2
}\|x_0\|_{\alpha}$ in the second inequality.
Choose $$\gamma\geq\gamma(2(1-\beta_0)p_2)\vee 2p_3\vee\gamma(2) \vee \gamma(4)\vee4\vee\gamma(2(p_0-1))$$
for some $p_0$ satisfying $p_0>-\frac{2}{\alpha}$, which will be used later.
Combining the above estimates we obtain that on $E_{K,\gamma}$
$$\aligned\int_1^t\|Y_1Z\|_{B^\alpha_{4,\infty}}^{2}ds\leq&C(p_1,p_2, \|x_0\|_{\alpha},K)(1+t) .\endaligned$$
Similarly by Lemma 3.5 we have for even $p>1$ with $\gamma\geq\gamma(p)$ that on $E_{K,\gamma}$
$$\aligned&\int_1^t\|\tilde{Y}_1(s)\|_{L^p}^pds\\\leq& C(p)[\int_1^t\|\tilde{Y}(s)\|_{L^p}^pds+\int_1^t\|e^{sA}x_1\|_{L^p}^pds]\\\leq& C(p,\|x_0\|_{\alpha},\|x_1\|_{\alpha})\int_1^t [1+\|Z\|_{\alpha}^{\gamma(p)}+\|:Z^2:\|_{\alpha}^{\gamma(p)}+\|:Z^3:\|_{\alpha}^{\gamma(p)}]ds
+C(\|x_1\|_{\alpha})(1+t)\\&+C(p,\lambda,\|\underline{Z}\|_{\mathfrak{L}_1},\|x_0\|_\alpha,\|x_1\|_\alpha)\\\leq &C(p, \|x_0\|_{\alpha},\|x_1\|_{\alpha},K)t+C(p,\lambda,K,\|x_0\|_\alpha,\|x_1\|_\alpha),\endaligned$$
and
$$\aligned&\int_1^t\|\nabla \tilde{Y}_1(s)\|_{L^2}^{2\beta_0 p_1}ds\\\leq& C(p_1,\|x_0\|_{\alpha},\|x_1\|_{\alpha})\lambda\int_1^t [1+\|Z\|_{\alpha}^{\gamma(2)}+\|:Z^2:\|_{\alpha}^{\gamma(2)}+\|:Z^3:\|_{\alpha}^{\gamma(2)}]ds
+C(\|x_1\|_{\alpha})(1+t)\\&+C(p_1,\lambda,\|\underline{Z}\|_{\mathfrak{L}_1},\|x_0\|_\alpha,\|x_1\|_\alpha)\\\leq &C(p_1, \|x_0\|_{\alpha},\|x_1\|_{\alpha},K)t\lambda+C(p_1,\lambda,K,\|x_0\|_\alpha,\|x_1\|_\alpha).\endaligned$$
Then we have on $E_{K,\gamma}$,
$$\aligned\|u(t)\|_{L^2}^2\leq& \|u(1)\|_{L^2}^2\exp [\int_1^t 2(-\lambda+L(s))ds]\\\leq&\|u(1)\|_{L^2}^2\exp [-\lambda t+C(p_1,p_2,\|x_0\|_{\alpha}, \|x_1\|_{\alpha},K)t\lambda^{\frac{1}{p_1}}+C(p_1,p_2,\lambda,\|x_0\|_{\alpha}, \|x_1\|_{\alpha},K)].\endaligned$$
By (3.7) and (3.12) we have $\|u(1)\|_{L^2}^2\leq C(\lambda,\|x_0\|_{\alpha}, \|x_1\|_{\alpha},K)$ on $E_{K,\gamma}$. Then we choose $\lambda$ large enough so that there exist constants $C_0, C_1>0$ such that
$$\|u(t)\|_{L^2}^2\leq C_0e^{-C_1t}\rightarrow0 \textrm{ on } E_{K,\gamma}, \textrm{ as } t\rightarrow\infty.$$
On the other hand by  Lemmas 3.4 and 3.5  we know that for $p_0>-\frac{2}{\alpha}$ on $E_{K,\gamma}$ and $t>1$
$$\|Y_1^{2p_0-2}(t)\|_{L^1}+\|\tilde{Y}_1^{2p_0-2}(t)\|_{L^1}\leq C(p_0, \|x_0\|_{\alpha},\|x_1\|_{\alpha},K,\lambda)(1+t),$$
which by H\"{o}lder's inequality implies that
$$\|u(t)\|_{L^{p_0}}\leq \|u(t)\|_{L^2}(\|Y_1^{2p_0-2}(t)\|_{L^1}^{\frac{1}{2}}+\|\tilde{Y}_1^{2p_0-2}(t)\|_{L^1}^{\frac{1}{2}})\rightarrow0 \textrm{ on } E_{K,\gamma}, \textrm{ as } t\rightarrow\infty,$$
Thus Lemma 2.1 yields that
$$\|u(t)\|_{\alpha}\rightarrow0 \textrm{ on } E_{K,\gamma}, \textrm{ as } t\rightarrow\infty.$$

From the above we also obtain that for fixed $K>0$, there exists $R>0$ such that $\tau_R=\infty$ on $E_{K,\gamma}$.
It follows that
$$\Gamma_{x_0,x_1}(D)\geq P(E_{K,\gamma})>0.$$
Now by Theorem 4.1 the first result of Theorem 1.1 follows.

\textbf{Proof of weak convergence}

In the following we prove that for fixed $x\in \mathcal{C}^\alpha$, $P_t(x,dy)$ converges to $\nu$ weakly, where $P_t(x,dy)$ denote the distribution of the $X(t)$ starting from $x$. We use similar arguments as the proof of  [KS16, Theorem 2.7].

By similar argument as the proof of [RZZ15, Theorem 3.10] we have that the solution $X$ to the equation (3.2) is continuous with respect to initial value in $\mathcal{C}^\alpha$, which implies  the Feller property of the semigroup easily.

  Now for $x\in\mathcal{C}^\alpha$ we prove the tightness of $\{P_n(x,dy),n\geq 1\}$.  By Lemma 3.6 for every $\varepsilon>0$, $y\in \mathcal{C}^\alpha$,  $\alpha<\alpha'<-\frac{2}{p_0}$ for $p_0$ above, we can find a generalized  coupling $\Gamma^\varepsilon_{x,y}\in \tilde{\mathcal{C}}(\delta_{x}P_\infty,\delta_y P_\infty)$ as above such that  $$\Gamma^\varepsilon_{x,y}(D)\geq 1-\varepsilon/2,\quad \Gamma^\varepsilon_{x,y}(\lim_{n\rightarrow\infty}\|x_n-y_n\|_{\alpha'}=0)\geq 1-\varepsilon/2,$$ and $\Pi_1^\sharp\Gamma^\varepsilon_{x,y}=\delta_{x}P_\infty$, where $\delta_{x}P_\infty$ denote the law of the sequence $\{X(n)\}$ on $\mathcal{C}^\alpha_\infty$ starting from $x$. In fact, $\Gamma^\varepsilon_{x,y}$ is the law of $(X(x),\hat{X}(y))$ as before and we choose $E_{\gamma,K(\varepsilon)}$  such that $P( E_{\gamma,K(\varepsilon)})\geq 1-\varepsilon/2$ and $\lambda$ depends on $K(\varepsilon)$, which makes the coupling dependent on $\varepsilon$.

   Define a measure on $\mathcal{C}^\alpha_\infty\times \mathcal{C}^\alpha_\infty$
 $$\Gamma^\varepsilon(A)=\int \Gamma^\varepsilon_{x,y}(A)\nu(dy), \quad A\in \mathcal{M}(\mathcal{C}^\alpha_\infty)\times \mathcal{M}(\mathcal{C}^\alpha_\infty).$$
 We have
 \begin{equation}\Gamma^\varepsilon(D)\geq 1-\varepsilon/2,\quad \Gamma^\varepsilon(\lim_{n\rightarrow\infty}\|x_n-y_n\|_{\alpha'}=0)\geq 1-\varepsilon/2.\end{equation}
 Since $\Pi_1^\sharp\Gamma^\varepsilon=\delta_{x}P_\infty$, $\Pi_2^\sharp\Gamma^\varepsilon\ll \nu P_\infty$, we deduce that
 $\Gamma^\varepsilon \in \tilde{\mathcal{C}}(\delta_{x}P_\infty,\nu P_\infty)$.  Moreover we have  for  $\varepsilon>0$ there exists $\delta>0$
 such that \begin{equation}\Gamma^\varepsilon (y_n\in K_0^c)\leq \varepsilon/3, \quad n\geq1,\end{equation}
 if a compact set $K_0\subset \mathcal{C}^\alpha$ is chosen such that $\nu(K_0)\geq 1-\delta.$

 Since the embedding $\mathcal{C}^{\alpha'}\subset \mathcal{C}^{\alpha}$ is compact and by [RZZ15, Lemma 3.1] we have that for every $k\in\mathbb{N}$, $\int \|\phi\|_{\mathcal{C}^{\alpha'}}^k\nu(d\phi)<\infty$,
 we can choose $C$ large enough such that for compact sets $K_1:=\{\|\cdot\|_{\alpha'}\leq C\}$, $K_2:=\{\|\cdot\|_{\alpha'}\leq C+1\}$
  $$\nu(K_2)\geq\nu(K_1)\geq1-\delta.$$
 By (4.10) we know that there exists $D_1$ such that $\Gamma^\varepsilon(D_1^c)\leq \frac{2\varepsilon}{3}$
 and $\|x_n-y_n\|_{\alpha'}$ converges to $0$ uniformly on $D_1$. For $n$ large enough, we have
 $$\Gamma^\varepsilon(x_n\in K_2^c)\leq \Gamma^\varepsilon(\{x_n\in K_2^c\}\cap D_1)+\Gamma^\varepsilon(D_1^c)\leq \Gamma^\varepsilon(y_n\in K_1^c)+\Gamma^\varepsilon(D_1^c)\leq \varepsilon,$$
 where we used (4.10), (4.11) in the last inequality. Since $\Pi_1^\sharp\Gamma^\varepsilon =\delta_{x}P_\infty$ we deduce the tightness of $\{P_n(x,\cdot),n\geq1\}$.

 In the following we prove the weak convergence: If we assume that  $P_n(x,\cdot)$ does not weakly converge to $\nu$, there exists some probability measure $\nu_0\neq \nu$ and a subsequence $P_{m_k}(x,\cdot)$ converges to $\nu_0$ weakly. Fix a bounded Lipschitz continuous function $f:\mathcal{C}^\alpha\rightarrow\mathbb{R}$ such that $\int fd\nu_0\neq \int fd\nu$ and set $U_n=\frac{1}{n}\sum_{k=1}^n f(x_{m_k})$.
 Now we want to prove that  $U_n$ converges to $\int fd\nu$ in probability with respect to $\delta_xP_\infty$.
For every $\varepsilon>0$ as above and construct corresponding $\Gamma^\varepsilon $ such that (4.10) holds.

 By [KS16, Corollary 2.6] we have that $U_n$ converges to $\int fd\nu$ in probability with respect to $\nu P_\infty$, which implies that  $U_n$ converges to $\int fd\nu$ in probability with respect to $\Pi_2^\sharp\Gamma^\varepsilon $. In fact,
 for every subsequence $\{n_r\}$ there exists another subsequence $\{n_{r_l}\}$ such that $U_{n_{r_l}}$ converges to  $\int fd\nu$ $\nu P_\infty$-a.s.. Since $\Pi_2^\sharp\Gamma^\varepsilon \ll \nu P_\infty$, $U_{n_{r_l}}$ converges to  $\int fd\nu$ $\Pi_2^\sharp\Gamma^\varepsilon $-a.s..

  Since $f$ is bounded and Lipschitz continuous, by (4.10) we have
 \begin{equation}\Gamma^\varepsilon (\lim_{n\rightarrow\infty}|\frac{1}{n}\sum_{k=1}^nf(x_{m_k})-\frac{1}{n}\sum_{k=1}^nf(y_{m_k})|=0)\geq1-\varepsilon/2.\end{equation}
We have for every $\varepsilon_0>0$
 $$\aligned&\delta_xP_\infty(|U_n-\int fd\nu|<\varepsilon_0)\\=&\Gamma^\varepsilon(|U_n-\int fd\nu|<\varepsilon_0)\\\geq& 1-\Gamma^\varepsilon  (|\frac{1}{n}\sum_{k=1}^nf(x_{m_k})-\frac{1}{n}\sum_{k=1}^nf(y_{m_k})|+|\frac{1}{n}\sum_{k=1}^nf(y_{m_k})-\int fd\nu|\geq\varepsilon_0)\\\geq&1-\Gamma^\varepsilon  (|\frac{1}{n}\sum_{k=1}^nf(x_{m_k})-\frac{1}{n}\sum_{k=1}^nf(y_{m_k})|\geq\frac{\varepsilon_0}{2})-\Gamma^\varepsilon (|\frac{1}{n}\sum_{k=1}^nf(y_{m_k})-\int fd\nu|\geq\frac{\varepsilon_0}{2}),\endaligned$$
which combined with (4.12) and the fact that  $U_n$ converges to $\int fd\nu$ in probability with respect to $\Pi_2^\sharp\Gamma^\varepsilon $ implies that
$$\lim_{n\rightarrow\infty}\delta_xP_\infty(|U_n-\int fd\nu|<\varepsilon_0)\geq1-\varepsilon.$$
Since $\varepsilon$ is arbitrary we deduce that $U_n$ converges to $\int fd\nu$ in probability with respect to $\delta_xP_\infty$.

Moreover, $f$ is bounded, we have that $$\int U_n d\delta_xP_\infty\rightarrow\int fd\nu.$$

On the other hand $P_{m_k}(x,\cdot)$ converges to $\nu_0$ weakly, we have
$$\int U_n d\delta_xP_\infty=\frac{1}{n}\sum_{k=1}^n\int f(y)P_{m_k}(x,dy)\rightarrow\int fd\nu_0\neq \int fd\nu,$$
which is a contradiction finishing the proof of the second result.

\textbf{Proof of (4.1) and (4.2)}

In the following we only have to prove (4.1), (4.2).  Let $a_\varepsilon=a*\rho_\varepsilon$. By a similar argument as in the proof of [RZZ15, Lemma 3.4] we have for $p>1$
$$:\hat{Z}^2:=\lim_{\varepsilon\rightarrow0}  (\hat{Z}^2_\varepsilon-c_\varepsilon) \textrm{ in } L^p(\Omega,C((0,\infty);\mathcal{C}^\alpha),Q),$$
$$:Z^2:+2Za+a^2=\lim_{\varepsilon\rightarrow0}  ({Z}^2_\varepsilon+2Z_\varepsilon a_\varepsilon+a_\varepsilon^2-c_\varepsilon)\textrm{ in } L^p(\Omega,C((0,\infty);\mathcal{C}^\alpha),P).$$
Since $${Z}^2_\varepsilon+2Z_\varepsilon a_\varepsilon+a_\varepsilon^2-c_\varepsilon=\hat{Z}^2_\varepsilon-c_\varepsilon,$$
and $P\thicksim Q$,
we obtain that (4.1) holds $P$-a.s..
(4.2) can also be proved by taking the limit as $\varepsilon\rightarrow0$ for the following equation
$$ \hat{Z}^3_\varepsilon-3c_\varepsilon\hat{Z}_\varepsilon={Z}^3_\varepsilon+3Z_\varepsilon a_\varepsilon^2+3(Z_\varepsilon^2-c_\varepsilon) a_\varepsilon+a_\varepsilon^3-3c_\varepsilon Z_\varepsilon.$$

$\hfill\Box$
\vskip.10in

In the following we will prove Theorem 1.4. First we introduce a space $\mathcal{F}C_b^\infty$, which will be used in the proof of Theorem 1.4.

Let $E=H^{-1-\epsilon}_2, E^*=H^{1+\epsilon}_2$ for some $\epsilon>0$.  We denote their Borel $\sigma$-algebras by $\mathcal{B}(E), \mathcal{B}(E^*)$ respectively.
Define $$\mathcal{F}C_b^\infty=\{u:u(z)=f({ }_{E^*}\!\langle l_1,z\rangle_E,{ }_{E^*}\!\langle l_2,z\rangle_E,...,{ }_{E^*}\!\langle l_m,z\rangle_E),z\in E, l_1,l_2,...,l_m\in E^*, m\in \mathbb{N}, f\in C_b^\infty(\mathbb{R}^m)\},$$
and  for $u\in \mathcal{F}C_b^\infty$ and $l\in L^2(\mathbb{T}^2)$, $$\frac{\partial u}{\partial l}(z):=\frac{d}{ds}u(z+sl)|_{s=0},z\in E,$$
Let $Du$ denote the $L^2$-derivative of $u\in \mathcal{F}C_b^\infty$, i.e. the map from $E$ to $L^2(\mathbb{T}^2)$ such that $$\langle Du(z),l\rangle=\frac{\partial u}{\partial l}(z)\textrm{ for all } l\in L^2(\mathbb{T}^2), z\in E.$$

\vskip.10in
\no \emph{Proof of Theorem 1.4.} First we prove that $\nu$ satisfies (i) and (ii) in Theorem 1.4. (i) is obvious from [GlJ86, Sect. 8.6].
By [AR91, Theorem 7.11]  the logarithmic derivative of $\nu$ along $k$ is
$$\beta_k=2\langle z,Ak\rangle-2\langle a_1:z^3:-a_2z,k\rangle,$$
for $z\in E$, $k\in C^\infty(\mathbb{T}^2)$, which implies (ii) by using [AR90, Corollary 4.8].

Let $\nu_0$ be the measure  satisfying (i), (ii) in Theorem 1.4. From (ii) we calculate the logarithmic derivative of  $\nu_0$: For $u\in \mathcal{F}C_b^\infty$, $k\in C^\infty(\mathbb{T}^2)$
$$\aligned\int \frac{\partial u}{\partial k}d\nu_0=&\int \lim_{t\rightarrow0}\frac{u(z+tk)-u(z)}{t}d\nu_0\\=&\lim_{t\rightarrow0}\int\frac{u(z+tk)-u(z)}{t}d\nu_0
\\=&\lim_{t\rightarrow0}\int\frac{(a_{-tk}(z)-1)u(z)}{t}d\nu_0\\=&\int\lim_{t\rightarrow0}\frac{a_{-tk}(z)-1}{t}u(z)d\nu_0,\endaligned$$
where in the second  equality we used that $u\in \mathcal{F}C_b^\infty$ and the dominated convergence theorem, and in the last equality we used (i) and [GlJ86, Section 8.6] to deduce the uniform integrability of $a_{tk}$. This implies  the logarithmic derivative of $\nu_0$ is the same as that of $\nu$.
 Hence  by [AR91]   the diffusion process $X^{\nu_0}$ obtained from the Dirichlet form $\mathcal{E}^0_{\nu_0}$ also satisfies (1.1) and $\nu_0$ is an invariant measure for $X^{\nu_0}$. Here  $\mathcal{E}^0_{\nu_0}$ is the closure of the pre-Dirichlet form
$$\mathcal{E}_{\nu_0}(u,v):=\frac{1}{2}\int_E \langle Du, Dv\rangle_{L^2}d\nu_0,$$
defined for  $u,v\in\mathcal{F}C_b^\infty$ (see [AR91]).
Moreover, by (i) we know that Lemma 3.6 in [RZZ15] also holds for $\nu_0$.  Furthermore, the same argument as in the proof of [RZZ15, Theorems 3.9] implies that $X^{\nu_0}$ also satisfies the shifted equation (3.2). By the uniqueness of the solution to (3.2) (see Theorem 3.1), we know that $\nu_0$ is also an invariant measure for the solution to (3.2). By Theorem 1.1 the result follows.$\hfill\Box$

\vskip.10in
\no \emph{Proof of Theorem 1.5.} First we prove that $\nu$ satisfies (i) and (ii) in Theorem 1.4. (i) is obvious from [GlJ86, Sect. 8.6].
By [AR91, Theorem 7.11]  the logarithmic derivative of $\nu$ along $k$ is
$$\beta_k=2\langle z,Ak\rangle-2\langle a_1:z^3:-a_2z,k\rangle,$$
for $z\in E$, $k\in C^\infty(\mathbb{T}^2)$, which implies (ii) by direct calculations.

Let $\nu_0$ be the measure  satisfying (i), (ii) in Theorem 1.5. From (ii)  we calculate the logarithmic derivative of  $\nu_0$:
We follow the proof of [BR95, Theorem 3.10]:
By (ii) we have $\int Lud\nu_0=0$ for $u\in \mathcal{F}C_b^\infty$. Hence for all $u,v \in \mathcal{F}C_b^\infty$
$$0=\int L(uv)d\nu_0=2\int uLvd\nu_0+\int \langle Du, Dv\rangle_{L^2} d\nu_0,$$
i.e.,
\begin{equation}-\int u Lvd\nu_0=\frac{1}{2}\int\langle Du, Dv\rangle_{L^2} d\nu_0.\end{equation}
Let $g_n\in C_b^\infty(\mathbb{R})$, $n\in\mathbb{N}$, such that $g_n(t)=t$ on $[-n,n]$ and $\sup\{|g_n'(t)|+|g_n''(t)|:n\in\mathbb{N}, t\in\mathbb{R}\}<\infty$.
Let $k\in C^\infty(\mathbb{T}^2)$. Applying (4.10) to $v:=g_n(k)$ we can take $n\rightarrow\infty$ according to the dominated convergence theorem, and since
$$L(g_n(k))=g_n''(k)\|k\|^2_{L^2}+g_n'(k)(\langle z,Ak\rangle-\langle a_1:z^3:-a_2z,k\rangle),$$
we obtain that
$$\aligned\int \frac{\partial u}{\partial k}d\nu_0=&-\int\beta_kud\nu_0. \endaligned$$
 Then we can conclude that the  the logarithmic derivative of $\nu_0$ along $k$ is the same as that of $\nu$. Hence  by the same proof as that for Theorem 1.4, the result follows. $\hfill\Box$

\vskip.10in
In the following we only prove Corollary 1.7. Corollary 1.8 can be obtained similarly.
\vskip.10in

\no \emph{Proof of Corollary 1.7.} Assume that $\nu$ can be written as a convex combination of two
probability measures $\mu_1$ and $\mu_2$ in $\mathcal{M}^a$. Then $\mu_1$ and $\mu_2$ are
absolutely
continuous w.r.t. to $\nu$ with bounded densities and hence are also absolutely
continuous
w.r.t.  the Gaussian measure $\mu$ with $p$-integrable densities for some $p>1$. By Theorem 1.4  $\mu_1=\mu_2=\nu$. So, $\nu$ is extreme in the
set $\mathcal{M}^a$. $\hfill\Box$


\begin{thebibliography}{99}
\bibitem[AKR97]{}S. Albeverio, Y. G. Kondratiev, M. R\"{o}ckner, Ergodicity for the Stochastic Dynamics of Quasi-invariant Measures with Applications to Gibbs States, Journal of Functional Analysis, 1997, 149(2), 415-469
\bibitem[AR90]{} S. Albeverio and M. R\"{o}ckner, Classical Dirichlet forms on topological vector spaces–
closability and a Cameron-Martin formula, Journal of Functional Analysis. 88 395-436, (1990)
\bibitem[AR91]{}S. Albeverio, M. R\"{o}ckner, Stochastic differential equations in infinite
dimensions: Solutions via Dirichlet forms, Probab. Theory Related Field 89 347-386 (1991).
\bibitem[BCD11]{} H. Bahouri, J.-Y. Chemin, R. Danchin,  Fourier analysis and nonlinear
partial differential equations, vol. 343 of Grundlehren der Mathematischen
Wissenschaften [Fundamental Principles of Mathematical Sciences]. Springer, Heidelberg,
2011.
\bibitem[BR95]{} V.I. Bogachev, M. R\"{o}ckner,  Regularity of Invariant Measures on Finite and Infinite Dimensional Spaces and Applications, Journal of Functional Analysis, 133, 1, 168-223, (1995)
\bibitem[CC13]{}R\'{e}mi Catellier, Khalil Chouk, Paracontrolled Distributions and the 3-dimensional Stochastic Quantization Equation, arXiv:1310.6869
\bibitem[DD03]{}G. Da Prato, A. Debussche, Strong solutions to the stochastic quantization equations, Ann.
Probab., 31(4):1900-1916, (2003)
\bibitem[DZ96]{} G. Da Prato, J. Zabczyk, Ergodicity for Infinite Dimensional Systems, London Mathematical Society Lecture Notes, n. 229, Cambridge University Press (1996)
\bibitem[GIP13]{} M. Gubinelli, P. Imkeller, N. Perkowski, Paracontrolled distributions and singular PDEs,  Forum Math. Pi 3 no. 6(2015)
\bibitem[GlJ86]{}  J.  Glimm, A. Jaffe : Quantum physics: a functional integral point of view. New
York Heidelberg Berlin: Springer (1986)
\bibitem[GRS75]{}F. Guerra, J. Rosen, B. Simon: The $P(\Phi)_2$ Euclidean quantum field theory
as classical statistical mechanics. Ann. Math. 101, l 11-259 (1975)
\bibitem[Hai14]{} M. Hairer, A theory of regularity structures. Invent. Math. , 198(2):269-504, (2014).
\bibitem[HMS11]{} M. Hairer, J. C. Mattingly, M. Scheutzow, Asymptotic coupling and a general form
of Harris theorem with applications to stochastic delay equations, Probability Theory
and Related Fields 149, 1-2, 223-259, (2011)
\bibitem[JLM85]{} G. Jona-Lasinio and P. K. Mitter. On the stochastic quantization of field theory. Comm.
Math. Phys., 101(3):409-436, (1985).
\bibitem[KS16]{}A. Kulik,
M. Scheutzow, Generalized couplings and convergence of transition probabilities, http://arxiv.org/abs/1512.06359v2
\bibitem[Kurz07]{}T. G. Kurtz, The Yamada-Watanabe-Engelbert theorem for general stochastic equations and inequalities, {Electronic Journal of Probability}. {12} 951-965, (2007)
    \bibitem[MR92]{} Z. M. Ma, and M. R\"{o}ckner, ”Introduction to the theory of (non-symmetric) Dirichlet
forms,” Springer-Verlag, Berlin/Heidelberg/New York, (1992)
    \bibitem[MR99]{}R. Mikulevicius, B. Rozovskii, Martingale problems for stochasic PDE's. In Stochastic
partial differential equations: six perspectives, volume 64 of Math. Surveys Monogr.
243-325. Amer. Math. Soc., Providence, RI, (1999)
 \bibitem[MW15]{}J. Mourrat, H. Weber,  Global well-posedness of the dynamic   $\Phi^4$ model in the plane, arXiv:1501.06191v1, to appear in The Annals of Probability
 \bibitem[PW81]{} G. Parisi,  Y. S. Wu. Perturbation theory without gauge fixing. Sci. Sinica 24,
no. 4, (1981), 483–496.
 \bibitem[Re95]{} S. Resnick, Dynamical Problems in Non-linear Advective Partial Differential Equations, PhD thesis, University of Chicago, Chicago (1995)
    \bibitem[R86]{ }M. R\"{o}ckner, Specifications and Martin boundaries for $P(\phi)_2$-random fields.
\emph{ Commun. Math. Phys.}\textbf{106}, 105-135 (1986)
\bibitem[RZ92] { }M. R\"{o}ckner and T.S. Zhang, Uniqueness of Generalized Schr\"{o}dinger Operators
and Applications, {Journal of Functional Analysis}. \textbf{105} 187-231 (1992)
\bibitem[RZZ15] { }M. R\"{o}ckner, R. Zhu, X. Zhu, Restricted Markov unqiueness for the
stochastic quantization of $P(\phi)_2$ and its
applications,  arXiv:1511.08030 (2015)
\bibitem[RZZ15a] { }M. R\"{o}ckner, R. Zhu, X. Zhu, Sub- and supercritical stochastic quasi-geostrophic equation, The Annals of Probability
43,  3, 1202-1273, (2015)
\bibitem[S74]{}B. Simon, The $P(\phi)_2$ Euclidean (Quantum) field theory. Princeton: Princeton
University Press (1974)
\bibitem[S85]{}W. Sickel, Periodic spaces and relations to strong summability of multiple Fourier series. Math. Nachr.
124, 15-44 (1985)

\bibitem[SW71]{}E. M. Stein, G. L. Weiss, Introduction to Fourier Analysis on Euclidean Spaces,
Princeton University Press, (1971)
\bibitem[Tri78]{} H. Triebel,  Interpolation theory, function spaces, differential operators. North-Holland Mathematical
Library 18. North-Holland Publishing Co. Amsterdam-New York 1978.
\bibitem[Tri06]{}H. Triebel, Theory of function spaces III. Basel, Birkh\"{a}user, (2006)




\end{thebibliography}
\end{document}